\documentclass[11pt]{article}
\usepackage{fullpage}
\usepackage{euscript}
\usepackage{amssymb,amsmath,exscale,relsize}
\usepackage{amsfonts}
\usepackage{amsbsy}
\usepackage{amsmath}
\usepackage{epsfig}
\usepackage{amsthm}
\usepackage{amscd}
\usepackage{amstext}
\usepackage{graphicx,dsfont,multirow,ctable,colortbl,rotating,slashbox,float}

\usepackage[pdftex]{hyperref}
\makeatletter
\renewcommand\section{\@startsection {section}{1}{\z@}
                                   {-3.5ex \@plus -1ex \@minus -.2ex}
                                   {2.3ex \@plus.2ex}
                                   {\normalfont\large\bfseries}}

\renewcommand\subsection{\@startsection{subsection}{2}{\z@}
                                     {-3.25ex\@plus -1ex \@minus -.2ex}
                                     {1.5ex \@plus .2ex}
                                     {\normalfont\bfseries}}

\newlength{\apb@width}
\newcommand{\autoparbox}[2][c]{\settowidth{\apb@width}{#2}\parbox[#1]{\apb@width}{#2}}

\makeatother

\newcommand{\bea}{\begin{eqnarray}}
\newcommand{\eea}{\end{eqnarray}}
\newcommand{\bee}{\begin{eqnarray*}}
\newcommand{\eee}{\end{eqnarray*}}
\newcommand{\al}{\begin{align*}}
\newcommand{\eal}{\end{align*}}
\newcommand{\be}{\begin{equation}}
\newcommand{\ee}{\end{equation}}
\newcommand{\eq}[1]{(\ref{#1})}
\newcommand{\bem}{\begin{pmatrix}}
\newcommand{\eem}{\end{pmatrix}}

\def\a{\alpha}

\def\c{\gamma}
\def\d{\delta}

\def\e{\epsilon}   
\def\f{\phi}               
  
\def\g{\gamma}
\def\h{\eta}

\def\inf{\infty}

\def\l{\lambda}
\def\m{\mu}

\def\p{\pi}   
\def\pa{\partial}       
\def\r{\rho}                                     
\def\s{\sigma}                                   
\def\t{\tau}
\def\th{\theta}
\def\til{\tilde}

\def\D{\Delta}
\def\F{\Phi}
\def\G{\Gamma}

\def\L{\Lambda}
\def\O{\Omega}

\def\U{\Upsilon}

\def\jac{\operatorname{jac}}
\def\reg{\operatorname{reg}}

\def\ex{\operatorname{e}}
\def\EG{\operatorname{{\bf EG}}}
\def\tr{\operatorname{tr}}
\def\tpi{2 \pi i}
\def\SL{\operatorname{\textsl{SL}}}
\def\Sp{\operatorname{\textsl{Sp}}}
\def\Id{\rm Id}
\def\tr{\operatorname{tr}}

\def \H {{\mathbb H}}

\def \Z {{\mathbb Z}}
\def \C {{\mathbb C}}
\def \R {{\mathbb R}}

\def\PSL{\operatorname{\textsl{PSL}}}
\def\Or{\operatorname{\textsl{O}}}
\def\Sp{\operatorname{\textsl{Sp}}}
\def\U{\operatorname{\textsl{U}}}
\def\SU{\operatorname{\textsl{SU}}}

\newtheorem{thm}{Theorem}[section]

\newtheorem{prop}[thm]{Proposition}
\newtheorem{conj}[thm]{Conjecture}

\theoremstyle{remark}

\begin{document}
\begin{center}

\vspace{1cm} 
{\fontsize{19}{0}\selectfont Ê\bf The Largest Mathieu Group\\  \vspace{ .3cm} and (Mock) Automorphic Forms\footnote{Submitted to the AMS Proceedings of Symposia in Pure Mathematics (String-Math 2011).}}
\vspace{1.5cm}

Miranda C. N. Cheng$~^{\flat, \natural, \sharp }$ and John F. R. Duncan$~ ^\dagger$

\vspace{0.8cm}

{\it 
$^\flat$Department of Mathematics, Harvard University,
\\
Cambridge, MA 02138, USA \\
$^\natural$Department of Physics, Harvard University,
\\
Cambridge, MA 02138, USA \\
$^\sharp$Universit\'e Paris 7, UMR CNRS 7586, Paris, France\\
$^\dagger$ Department of Mathematics, Case Western Reserve University,\\
Cleveland, OH 44106, USA
}

\vspace{0.6cm}

\end{center}

\begin{abstract}
We review the relationship between the largest Mathieu group and various modular objects, including recent progress on the relation to mock modular forms.  We also review the connections between these mathematical structures and string theory on $K3$ surfaces.   
 \end{abstract}

\setcounter{page}{1}
\tableofcontents

\section{Cusp Forms and $M_{24}$} 
\label{Cusp}
\setcounter{equation}{0}

We start by recalling some facts about the largest Mathieu group, $M_{24}$. The group $M_{24}$ may be characterised as the automorphism group of the unique doubly even self-dual binary code of length $24$ with no words of weight $4$, also known as the {\em (extended) binary Golay code}. 
In other words, there is a unique (up to permutation) set, ${\mathcal G}$ say, of length $24$ binary codewords (sequences of $0$'s and $1$'s) such that any other length $24$ codeword has even overlap with all the codewords of ${\mathcal G}$ if and only if this word itself is in ${\mathcal G}$, and the number of $1$'s in each codeword of ${\mathcal G}$ is divisible by but not equal to 4. The group of permutations of the $24$ coordinates that preserves the set ${\mathcal G}$ is the sporadic group $M_{24}$.  See, for instance, \cite{sphere_packing}. 
 \begin{table}[h] \centering  \begin{tabular}{ccccccc}
 \toprule
$[g]$ & cycle shape & $\eta_g(\t)$ & $k_g$ &$n_g$& $N_g$ &$h_g$ \\\midrule
$1A$ & $1^{24}$ & $ \eta^{24}(\t) $ & 12 &1  &1 &1\\
$2A$ & $1^{8}2^8$ & $\eta^{8}(\t)\eta^{8}(2\t)$ & 8&2&2 &1\\
$2B $& $2^{12}$ & $\eta^{12}(2\t)$ & 6 &2& 4&2\\
$3A$&$1^63^6$&$\h^6(\t)\h^6(3\t)$&6&3&3&1\\
$3B$&$3^8$&$\h^8(3\t)$&4&3&9&3\\
$4A$& $2^4 4^4$ & $\h^4(2\t)\h^4(4\t)$ &4&4&8&2\\
$4B$&$1^4 2^2 4^4$&$\h^4(\t)\h^2(2\t)\h^4(4\t)$&5&4&4&1\\
$4C$&$4^6$&$\h^6(4\t)$&3&4&16&4\\
$5A$&$1^45^4$&$\h^4(\t)\h^4(5\t)$&4&5&5&1\\
$6A$&$1^22^23^26^2$&$\h^2(\t)\h^2(2\t)\h^2(3\t)\h^2(6\t)$&4&6&6&1\\
$6B$&$6^4$&$\h^4(6\t)$&2&6&36&6\\
$7AB$&$1^3 7^3$&$\h^3(\t)\h^3(7\t)$&3&7&7&1\\
$8A$&$1^2 2^14^18^2$&$\h^2(\t)\h(2\t)\h(4\t)\h^2(8\t)$&3&8&8&1\\
$10A$&$2^210^2$&$\h^2(2\t)\h^2(10\t)$&2&10&20&2\\
$11A$&$1^2 11^2$&$\h^2(\t)\h^2(11\t)$&2&11&11&1\\
$12A$&$2^14^16^112^1$&$\h(2\t)\h(4\t)\h(6\t)\h(12\t)$&2&12&24&2\\
$12B$&$12^2$&$\h^2(12\t)$&1&12&144&12\\
$14AB$&$1^1 2^1 7^114^1$&$\h(\t)\h(2\t)\h(7\t)\h(14\t)$&2&14&14&1\\
$15AB$&$1^1 3^1 5^115^1$&$\h(\t)\h(3\t)\h(5\t)\h(15\t)$&2&15&15&1\\
$21AB$&$3^1 21^1$&$\h(3\t)\h(21\t)$&1&21&63&3\\
$23AB$&$1^1 23^1$&$\h(\t)\h(23\t)$&1&23&23&1\\
 \bottomrule
  \end{tabular}
   \caption{\label{examples_eta} \footnotesize{The cycle shapes, weights $(k_g)$, levels $(N_g)$ and orders $(n_g)$ of the 26 conjugacy classes of the sporadic group $M_{24}$. The length of the shortest cycle is $h_g=N_g/n_g$. The naming of the conjugacy classes follows  the ATLAS convention \cite{atlas}. We write $7AB$, for example, to indicate that the entries of the incident row are valid for both the conjugacy classes $7A$ and $7B$.
}}
  \end{table}

As such, $M_{24}$ naturally admits a permutation representation of degree $24$, which we denote $R$ and call the {\em defining representation} of $M_{24}$, and via $R$ we may assign a {\em cycle shape} to each of its elements. For example, to the identity element we associate the cycle shape $1^{24}$, to an element of $M_{24}$ that is a product of $12$ mutually commuting transpositions we associate the cycle shape $2^{12}$, and so on. More generally, any cycle shape arising from an element of $M_{24}$ (or $S_{24}$, for that matter) is of the form
$$
{i_1}^{\ell_1} {i_2}^{\ell_2} \dotsi{i_r}^{\ell_r} ,\quad \sum_{s=1}^{r} \ell_{s} \,i_{s} =24\;,
$$
for some $\ell_s\in\Z^+$ and $1\leq i_1<\cdots<i_r\leq23$ with $r\geq 1$. 
Clearly, the cycle shape of an element $g\in M_{24}$ depends only on its conjugacy class, denoted by $[g]$, although different conjugacy classes can share the same cycle shape.
We write
\vspace{-0.1cm}
$$
\chi(g) = \tr_R g
$$
for the character attached to this $24$-dimensional representation of $M_{24}$. The decomposition of $R$ into irreducible representations of $M_{24}$ is given by  
$$
 R = \rho_1 \oplus \rho_{23} \;
$$
(cf. Table \ref{M24}). The value of the character $\chi(g)$ equals the number of fixed points of $g$ in the action on the set of 24 points. In particular, note that $\chi(g)=\ell_1$ in case $i_1=1$ and $\chi(g)=0$ otherwise.

It turns out that the cycle shapes of $M_{24}$ have many special properties that will be important for the understanding of the modular properties of the associated McKay--Thompson series which we will discuss shortly. 
First, the $M_{24}$ cycle shapes are privileged in that they are all of the so-called {\em balanced type} \cite{conway_norton}, meaning that for each $g\in M_{24}$ there exists a positive integer $N_g$ such 
that if $\prod i_s^{\ell_s}$ is the cycle shape associated to $g$ then
$$
	\prod_si_s^{\ell_s}=\prod_s\left(\frac{N_g}{i_s}\right)^{\ell_s}\;.
$$
We will refer to the number $N_g$ as the {\em level} of the $g$. 

If $g$ has cycle shape $i_1^{\ell_1}\cdots i_r^{\ell_r}$ then the order of $g$ is the least common multiple of the $i_s$'s. 
A second special property of $[g]\subset M_{24}$ is that the order of $g$ coincides with the length $i_r$ of the longest cycle in the cycle shape. 
Henceforth we will denote $n_g=i_r$. 

Finally, observe that for all $g\in M_{24}$ the level $N_g$ defined above equals the product of the shortest and the longest cycle. 
Hence we have $n_gh_g  = N_g$ where $h_g$ denotes the length of the shortest cycle in the cycle shape. 
Moreover, we also have the property $h_g \lvert n_g$ and $h_g\lvert 12$. This is reminiscent of the fact, significant in monstrous moonshine \cite{conway_norton}, that the normaliser of $\Gamma_0(N)/\langle\pm \Id\rangle$ in $\PSL_2(\R)$ (cf. (\ref{congruence1})) is a group containing a conjugate of $\Gamma_0(n/h)/\langle\pm \Id\rangle$ where $n=N/h$ and $h$ is the largest divisor of $24$ such that $h^2|N$. We also set $k_g=\sum_{s=1}^r \ell_s$/2 to be half of the total number of cycles and call it the {\em weight} of $g$. 
Of course, $N_g$, $n_g$, $h_g$ and $k_g$ depend only on the conjugacy class $[g]$ containing $g$. All the values of these parameters can be found in Table \ref{examples_eta}.

To each element $g\in M_{24}$ we can attach an {\em eta-product}, to be denoted $\eta_g$, which is the function on the upper half-plane given by
\begin{gather}
	\eta_g(\tau)=\prod_s\eta(i_s\tau)^{\ell_s}
\end{gather}
where $\prod_si_s^{\ell_s}$ is the cycle shape attached to $g$, and $\eta(\tau)$ is the Dedekind eta function satisfying $\eta(\tau)=q^{1/24}\prod_{n\in\Z^+}(1-q^n)$ for $q=e(\tau)$, where here and everywhere else in the paper we use the shorthand notation 
$$e(x)=e^{2\p i x}\;.$$ 
As was observed in \cite{Mason,DummitKisilevskyMcKay}, the eta-product $\eta_g$ associated to an element $g\in M_{24}$ (or rather, to its conjugacy class $[g]$) is a cusp form of  weight $k_g$ for the group $\G_0(N_g)$, with a {\em Dirichlet character} $\varsigma_g$ that is trivial if the weight $k_g$ is even and is otherwise defined, in terms of the Jacobi symbol $(\frac{n}{m})$, by
$$
\varsigma_g(\g)=
\begin{cases}
	\left(\frac{N_g}{d}\right) (-1)^{\frac{d-1}{2}},& 
	\text{$d$ odd,} \\
	\left(\frac{N_g}{d}\right), & 
	\text{$d$ even,} 
\end{cases}
$$
in case $d$ is the lower right entry of $\g \in \G_0(N_g)$. 
Let's recall that 
\be\label{congruence1}
\G_0(N)=\left\{\g \,\Big\lvert\; \g = \bem a&b \\ c&d \eem \in \SL_2(\Z)\,,\; c\equiv  0 \;\;{\rm mod}\;\; N \right\}\;,
\ee
and a cusp form is a modular form which vanishes in the limit as $\t \to i\inf$. 

Interestingly, these cusp forms point to the existence of an infinite-dimensional $\Z$-graded module for $M_{24}$; a fact we will now explain. It is by now a routine task to construct a module for $M_{24}$ whose graded dimension is the reciprocal of the cusp form $\eta^{24}(\t)$. 
To see this, let us consider the defining 24-dimensional representation $R$. To obtain a $q$-series we also need an extra $\Z$-grading. Hence we will take a copy of this representation, denoted by $R_k$, for each positive integer $k$ and consider the direct sum
$$
V= R_{1} \oplus R_{2} \oplus  \dotsi  = \bigoplus_{k=1}^\inf R_{k}\;.
$$
 Furthermore we would like to consider the graded Fock space ${\cal H}$ built from the above infinite-dimensional vector space $V$
 $$
 {\cal H} = \bigoplus_{n=0}^\inf S^n V = \bigoplus_{k=0}^\inf {\cal H}_{k} \quad,\quad  S^nV= V^{\otimes n}/S_n \;,$$
where the $\Z$-grading is naturally inherited from that of $V$. 
The associated partition function is then given by 
$$
Z(\t) = \sum_{k=0}^\inf q^{k-1} \text{dim} {\cal H}_{k} = \frac{1}{\eta^{24}(\t)}\;.
$$

In the language of vertex operator algebras, this is simply the partition function of the conformal field theory ({CFT}) with 24 free chiral bosons. Namely,  we consider 24 pairs of creation and annihilation operators $a_{-k}^{(i)}, a_{k}^{(i)}, i=1,\dotsi,24$ for each positive integer $k$, which furnish the 24-dimensional representation $R_k$ and satisfy the commutation relation
$$
[a_k^{(i)}, a_\ell^{(j)}] = \d_{ij} \d_{k+\ell,0}\;. 
$$
The quantum states of this theory are built from the unique vacuum state $\lvert 0\rangle$, characterised by the condition $a_n^{(i)}\lvert 0\rangle = 0 $ for all $n\geq 0, i=1,\dotsi,24$, and take the form
$$
a_{-n_1}^{(i_1)}\dotsi a_{-n_N}^{(i_N)} \lvert 0\rangle \quad,\quad n_k >0 , \;\; 1 \leq i_k \leq 24 \;,
$$
while the Hamiltonian of the theory is given by the operator
$$
\hat H= \sum_{i=1}^{24}\sum_{n=1}^\inf n\, a_{-n}^{(i)}\, a_n^{(i)} -1\;.
$$
Hence the partition function of this theory is indeed given by the modular form
$$
Z(\t) =  \tr_{\cal H} q^{\hat H} = \frac{1}{\eta^{24}(\t)} \;.
$$
Notice that, physically the $(-1)$-shift in the Hamiltonian operator corresponds to the ground state energy in the quantum theory with ${\cal H}$ as the space of quantum states (the {\em Hilbert space}).

From the above construction, in particular the fact that the $M_{24}$-action  commutes with the $\Z$-grading,
 we conclude that not only ${\cal H}$ but also each of its components 
${\cal H}_{k}$ are  $M_{24}$-modules. Hence, for each element $g\in M_{24}$ we can consider the following {\em twisted} partition function, or {\em McKay-Thompson series}
\be\label{mckay_thompson_eta1}
Z_g(\t)  =  \tr_{\cal H} g q^{\hat H}=\sum_{k=0}^\inf q^{k-1} \big( \tr_{{\cal H}_{k}} g\big) \;.
\ee
To be more precise, in the above formula $g$ really stands for the corresponding linear map acting on the representation in question. 
Note that taking the identity element $g={\mathds{1}}$, we recover the usual partition function $ \tr_{\cal H}  q^{\hat H} = Z(\t)$ encoding the graded dimension of the representation.  Moreover, it is clear from their definition that these  McKay-Thompson series depend only on the conjugacy class $[g]$ of the group element. Namely, we have 
$$
Z_g(\t) = Z_{hgh^{-1}}(\t)\;.
$$

From the description of ${\cal H}$, it is now easy to prove that they are given by the eta-products as
\be\label{twisting1}
Z_g(\t)  = q^{-1} \exp\left(\sum_{k=1}^\inf\tr_V   \frac{(q^{\hat H} g)^k}{k}  \right) = q^{-1} \exp\left( \sum_{n=1}^\inf  \sum_{k=1}^\inf \frac{q^{nk} \,\tr_R(g^k)}{k}  \right)  = \frac{1}{\eta_g(\t)}\;.
\ee
In the path integral language, the summing over $k$ in the expression for $\log Z_k$ can be thought of as a sum over how many times the time circle of the world-volume wraps around the time circle in the target space.

In other words, we can equate the McKay-Thompson series of the infinite-dimensional $M_{24}$-module ${\cal H}$, given by the twisted partition function of the physical theory $Z_g(\t)$, to a weight $-k_g$, level $N_g$ modular form $1/\eta_g(\t)$. 
The first few Fourier coefficients of the twisted partition functions and the corresponding $M_{24}$-modules ${\cal H}_k$ in terms of its decomposition into irreducible representations can be found in Table \ref{eta_prod_decomp}.

A few comments are in order here. First, historically the relation between the $\h_g(\t)$ and $M_{24}$ was phrased in terms of the eta-products themselves as opposed to their reciprocals as we described above. In that case, underlying the Fourier coefficients is an infinite-dimensional, $\Z$-graded {\it super} (or {\em virtual}) module of $M_{24}$, corresponding to considering instead the following $\Z/2$-graded Hilbert space
$$
\til {\cal H}=\til {\cal H}_{\bar{0}}\oplus \til {\cal H}_{\bar{1}},\quad
\til {\cal H}_{\bar{k}} =  \bigoplus_{n=0}^\inf  \L^{2n+k} V 
$$
which is the exterior (super)algebra of the infinite-dimensional vector space $V$, where $V$ is regarded as having odd parity \cite{Mason,RamanunjanEta}. 
A super $G$-module is just a $G$-module $W$ with a $G$-invariant $\Z/2$-grading $W=W_{\bar{0}}\oplus W_{\bar{1}}$. Given such a module we may consider the operator $F:W\to W$ for which the subspaces $W_{\bar{k}}$ are eigenspaces with eigenvalue $k$ for $F$. Then the super trace of an operator $A$ on $W$ is defined by setting
\[
{\rm str}_WA ={\rm tr}_W(-1)^FA.
\]
In physical terms such a $\Z/2$-grading often indicates the presence of fermions, and the operator $F$ is called the {\em fermion number operator}.

Second, as the alert reader may have noticed, starting from $\eta^{24}(\t)$ we can build not only an infinite-dimensional module ${\cal H}$ for $M_{24}$ but this same space (or vertex operator algebra) admits an action by the full symmetric group $S_{24}$, or even $\Or_{24}(\C)$\footnote{Although in general the corresponding twisted partition function won't be multiplicative if $g \in S_{24}$ but $g\notin M_{24}$. Of the 1575 possible partitions of $24$ only $30$ correspond to multiplicative eta-products \cite{DummitKisilevskyMcKay}. All of the 30 multiplicative eta-products are recovered from a permutation action of a group of the form $2^{24} M_{24}$ and 21 of them arise from the action of $M_{24}$; these are the 21 eta-products listed in Table \ref{examples_eta}. }. In this sense, the relation between these 21 eta-products and $M_{24}$ is not as distinguished as one might have hoped from the point of view of moonshine. Nevertheless, as we will see later in \S\ref{Siegel}, together with various weak Jacobi forms, mock modular forms and Siegel modular forms, they do tie together and form a closely-knitted net of modular objects attached to the sporadic group $M_{24}$.

Third, while we have seen that the eta-products $\eta_g(\t)$ are cusp forms of level $N_g$, from the following standard CFT arguments we expect a further, more refined modular property. To explain this, let us first discuss a certain natural generalisation of the above discussion:
Given a vertex operator algebra $\cal H$ and two commuting automorphisms $g$ and $h$ we expect there to exist a so-called {\em $h$-twisted module} ${\cal H}^{(h)}$ for $\cal H$ and one can consider the twisted twisted sector partition function 
\begin{gather}\label{eqn:Zgh}
Z\left(\t;\mbox{\fontsize{9}{10}\selectfont $g\,$}\underset{h}{\Box}\right)
\end{gather} 
defined as the (rationally) graded trace of $g$ on ${\cal H}^{(h)}$, thereby generalizing the twisted partition function $Z_g(\t)$ which is recovered by taking $h$ to be the identity. The significance of the twisted partition functions (\ref{eqn:Zgh}) in the context of monstrous moonshine was first observed by S. Norton in \cite{generalized_moonshine} and dubbed {\em generalised moonshine}. It is conjectured in loc. cit. that any such function is either constant or a generator for the field of functions on some discrete subgroup of $\PSL_2(\R)$ in case $\cal H$ is the moonshine module vertex operator algebra and $g$ and $h$ are a commuting pair of elements of the Monster simple group. An analogue of generalised moonshine for commuting pairs of elements of $M_{24}$ is discussed in \cite{Mason_genus0}. 
 
Recall that the equivalence between the Hamiltonian and Lagrangian formulations of quantum field theory implies that the insertion of the elements $g$ and $h$ as in (\ref{eqn:Zgh}) has an interpretation as altering the boundary condition of the path integral over the maps from an elliptic curve world-sheet into the target space. Identifying the modular group $\SL_2(\Z)$ as the mapping class group acting on the elliptic curve world-sheet, we conclude that the following must hold
$$
Z\left(\t;\mbox{\fontsize{9}{10}\selectfont $g\,$}\underset{h}{\Box}\right)
= 
\varepsilon(\gamma,g,h)\, 
Z\left(\t|_\g;\mbox{\fontsize{9}{10}\selectfont $g^a \!h^b\!\!$}\underset{g^c h^d}{\Box}\right),\quad \g =\bem a&b\\ c&d \eem \in \SL_2(\Z),
$$
where $\varepsilon(\gamma,g,h)$ is a phase. 
Specializing the above equation to the case that $h$ is the identity we conclude that the McKay-Thompson series $Z_g$ must be invariant up to a phase under the group $\Gamma_0(n_g)$ where $n_g$ is the order of the group.

Indeed, the eta-product $\eta_g$ also defines a cusp form of weight $k_g$ on the larger (or equal) group $\G_0(n_g) \supseteq \G_0(N_g)$ if we allow for a slightly more subtle multiplier system. Note that this change in multiplier system only comes into play when $g$ is in one of the fixed-point-free classes of $M_{24}$ with $\chi(g)=0$, for otherwise $\G_0(n_g) = \G_0(N_g)$.
Our convention is to say that a function $\xi:\G\to \C^*$ is a {\em multiplier system for $\G$ of weight $w$} in case the identity
\be\label{defn:multsys}
	\xi(\g\sigma)
	\jac(\g\sigma,\t)^{w/2}
	=
	\xi(\g)
	\xi(\sigma)
	\jac(\g,\sigma\t)^{w/2}
	\jac(\sigma,\t)^{w/2}
\ee
holds for all $\g,\sigma\in\G$ and $\t\in\H$. Writing $\g= \big(\begin{smallmatrix}a&b\\c&d\end{smallmatrix}\big)$, we have in the above formula $\g\t=\frac{a\t+b}{c\t+d}$, and the Jacobian $$\jac(\g,\t)=(c\t+d)^{-2}\;.$$
In detail, for $\g\in \G$, we define the {\it slash operator $\lvert_{\xi,w}$ of weight $w$ associated with the multiplier $\xi$} by setting
\be\label{slash}
(f|_{\xi,w}\g)(\t)
	=
	\xi(\g)f(\g\t)\jac(\g,\t)^{\mathsmaller{w/2}}
\ee
so that a holomorphic function $f: \H\to \C$ is a modular form of weight $w$ and multiplier $\xi$ for $\G$ if and only if $(f|_{\xi,w}\g)(\t)=f(\t)$ for all $\g\in \G$. 
In the present case we have 
\be\label{eqn:G0nvaretag}
\Big(\frac{1}{\h_g}\Big\lvert_{\xi_g,k_g} \g\Big)(\t) = \frac{1}{\h_g(\t)},\quad \text{for all}\quad \g \in \G_0(n_g),
\ee
where $\xi_g(\g) = \r_{n_g|h_g}(\g) \varsigma_g(\g)$ with 
\be\label{rho_multiplier}
\rho_{n|h}(\g)=\ex(-\tfrac{1}{(\g\infty-\g 0)nh})=\ex(-cd/nh)\;.
\ee
Note that $\rho_{n_g|h_g}$ is actually a character on $\G_0(n_g)$ since we have that $xy\equiv 1 \pmod{h_g}$ implies $x\equiv y\pmod{h_g}$ by virtue of the fact that $h_g=N_g/n_g$ is a divisor of $24$ for every $g\in M_{24}$ \cite[\S3]{conway_norton}. 

\section{Mock Modular Forms and $M_{24}$} 
\label{Mock}
\setcounter{equation}{0}

In the previous section we have seen the relation between $M_{24}$ and $\eta$-quotients. Recently, there has been an unexpected observation relating $M_{24}$ and a weight $1/2$ mock modular form. Namely, consider the $q$-series 
$$
H(\t)= 2 q^{-\frac{1}{8}} \left( -1 + 45 q+ 231 q^2 +\dots \right) = q^{-\frac{1}{8}} \left(-2  + \sum_{n=1}^\inf t_nq^n\right).
$$
Then the observation made in \cite{Eguchi2010} is that the first few $t_n$'s read
$$
2\times45,\, 2\times231,\, 2\times770,\,2\times2277,\,2\times 5796,\,\ldots
$$
and the integers $45$, $231$, $770$, $2277$ and $5796$ are dimensions of irreducible representations of $M_{24}$. 

This function $H(\tau)$, defined in  \eq{second_def_H}, enjoys a special relationship with the group $\SL_2(\Z)$; namely, it is a {\em weakly holomorphic mock modular form of weight $1/2$} on $\SL_2(\Z)$ with {\em shadow} $24 \,\eta(\tau)^3$ \cite{Eguchi2008,Atish_Sameer_Don}. 
This means that if we define the {\em completion}  $\hat H(\tau)$ of the holomorphic function $H(\tau)$ by setting 
\be\label{completion_H}
\hat{H}(\t)=H(\t)+24\, (4{i})^{-1/2} \int_{-\bar \t}^{\infty}(z+\t)^{-1/2}\overline{\eta(-\bar z)^3}{\rm d}z,
\ee
then $\hat{H}(\t)$ transforms as a modular form of weight $1/2$ on $\SL_2(\Z)$ with multiplier system conjugate to that of $\eta(\tau)^3$. In other words, we have
$$
\big(\hat{H}(\t)\lvert_{\e^{\mathsmaller{-3}},\mathsmaller{1/2}}\mathlarger{\g}\big)(\t)= \epsilon(\gamma)^{\mathsmaller{-3}}\hat{H}(\g\t)
\jac(\g,\t)^{\mathsmaller{1/4}}=\hat{H}(\t)
$$
for $\gamma\in \SL_2(\Z)$, where $\epsilon: \SL_2(\Z)\to \C^*$ is the multiplier system for $\eta(\t)$ satisfying 
$$\big(\h\lvert_{\e,\mathsmaller{1/2}}\mathlarger{\g}\big)(\t)=\eta(\t)\;.$$
See \S\ref{Modular Forms} for an explicit description of $\e$.

More generally, a holomorphic function $h(\t)$ on $\H$ is called a {\em (weakly holomorphic) mock modular form of weight $w$} for a discrete group $\G$ (e.g. a congruence subgroup of $\SL_2(\R)$) if it has at most exponential growth as $\t\to\alpha$ for any $\alpha\in \mathbb{Q}$, and if there exists a holomorphic modular form $f(\t)$ of weight $2-w$ on $\G$ such that $\hat{h}(\t)$, given by
\be\label{def_mock}
\hat{h}(\t)=h(\t)+\left(4{i}\right)^{w-1}\int_{-\bar \t}^{\infty}(z+\t)^{-w}\overline{f(-\bar z)}{\rm d}z,
\ee
is a (non-holomorphic) modular form of weight $w$ for $\G$ for some multiplier system $\psi$ say. In this case the function $f$ is called the {\em shadow} of the mock modular form $h$ and $\psi$ is called the multiplier system of $h$. Evidently $\psi$ is the conjugate of the multiplier system of $f$. The completion $\hat{h}(\t)$ satisfies interesting differential equations. For instance, completions of mock modular forms were identified as Maass forms in \cite{BriOno_MockThetaConj} and this led to a solution of the longstanding Andrews--Dragonette conjecture (cf. loc. cit.). As was observed in \cite{Atish_Sameer_Don} we have the identity
$$
	2^{1-w}\pi\Im(\t)^w\frac{\partial\hat{h}(\t)}{\partial\bar{\t}}=-2\pi i\overline{f(\t)}
$$
when $f$ is the shadow of $h$.

Given the observation regarding the first few Fourier coefficients of $H(\t)$ indicated above one would like to conjecture that the entire set of values $t_{n}$ for $n\in \Z^+$ encode the graded dimension of a naturally defined $\Z$-graded $M_{24}$ module $K=\bigoplus_{n=1}^{\inf} K_n$ with  ${\rm  dim}\,K_n\!=\!t_n\;.$ Of course, this conjecture by itself is an empty statement, since all positive integers can be expressed as dimensions of representations of any group, since we may always consider trivial representations. However, the fact that the first few $t_n$ can be written so nicely in terms of irreducible representations suggests that the $K_n$ should generally be non-trivial, and given any particular guess for a $M_{24}$-module structure on the $K_n$ we can test its merit by considering the twists of this mock modular form $H$ obtained by replacing the identity element in $\dim K_n=\tr_{K_n}\mathds{1}$ with an element $g$ of the group $M_{24}$. We would then call the resulting $q$-series
\begin{gather}
q^{-\frac{1}{8}} \left(-2  + \sum_{n=1}^\inf \tr_{K_n}gq^n\right)
\end{gather}
the {\em McKay--Thompson series} attached to $g$. A non-trivial connection between mock modular forms and $M_{24}$ arises if all such McKay--Thompson series of the $M_{24}$-module $K$ display interesting (mock) modular properties. In fact, since a function with good modular properties is generally determined by the first few of its Fourier coefficients, it is easier in practice to guess the McKay--Thompson series than it is to guess the representations $K_n$. Not long after the original observation was announced in \cite{Eguchi2010} candidates for the McKay--Thompson series had been proposed for all conjugacy classes $[g]\subset M_{24}$ in \cite{Cheng2010_1,Gaberdiel2010,Gaberdiel2010a,Eguchi2010a}, and with functions $\tilde{T}_g(\t)$ defined as in Table \ref{h_g} the following result was established.
\begin{table}[h!] \centering  \begin{tabular}{ccc}
 \toprule
$[g]$ & $\chi(g)$ & $\tilde{T}_g(\tau)$\\\midrule
$1A$ & $24$ & $0 $\\
$2A$ & ${8}$ & $16\Lambda_2$\\
$2B $&0& $-24\Lambda_2+8\Lambda_4=2\eta(\tau)^8/\eta(2\tau)^4$\\
$3A$&6& $6\Lambda_3$\\
$3B$&0& $2\eta(\tau)^6/\eta(3\tau)^2$\\
$4A$&0& $4\Lambda_2-6\Lambda_4+2\Lambda_8=2\eta(2\tau)^8/\eta(4\tau)^4$\\
$4B$&4& $4(-\Lambda_2+\Lambda_4)$\\
$4C$&0&$ 2\eta(\tau)^4\eta(2\tau)^2/\eta(4\tau)^2$\\
$5A$ &4&$ 2\Lambda_5$\\
$6A$&2&$ 2(-\Lambda_2-\Lambda_3+\Lambda_6)$\\
$6B$&0&$ 2\eta(\tau)^2\eta(2\tau)^2\eta(3\tau)^2/\eta(6\tau)^2$\\
$7AB$&3&$ \Lambda_7$\\
$8A$&2&$ -\Lambda_4+\Lambda_8$\\
$10A$& 0&$ 2\eta(\tau)^3\eta(2\tau)\eta(5\tau)/\eta(10\tau)$\\
$11A$& 2&$ 2(\Lambda_{11}-11\eta(\tau)^2\eta(11\tau)^2)/5$\\
$12A$& 0&$ 2\eta(\tau)^3\eta(4\tau)^2\eta(6\tau)^3/\eta(2\tau)\eta(3\tau)\eta(12\tau)^2$\\
$12B$& 0&$ 2\eta(\tau)^4\eta(4\tau)\eta(6\tau)/\eta(2\tau)\eta(12\tau)$\\
$14AB$& 1&$ (-\Lambda_2-\Lambda_7+\Lambda_{14}-14\eta(\tau)\eta(2\tau)\eta(7\tau)\eta(14\tau))/3$\\
$15AB$& 1&$ (-\Lambda_3-\Lambda_5+\Lambda_{15}-15\eta(\tau)\eta(3\tau)\eta(5\tau)\eta(15\tau))/4$\\
$21AB$& 0&$ (7\eta(\tau)^3\eta(7\tau)^3/\eta(3\tau)\eta(21\tau)-\eta(\tau)^6/\eta(3\tau)^2)/3$\\
$23AB$& 1&$ (\Lambda_{23}-23f_{23,1}+23f_{23,2})/11$\\
 \bottomrule
  \end{tabular}
   \caption{\label{h_g} \footnotesize{In this table we collect the data that via equation \eq{h_g_explicit} and \eq{twisted_elliptic_genus2}  define the weight $1/2$ (mock) modular forms $H_g(\t)$ and the weak Jacobi forms ${\cal Z}_g(\t,z)$. For $N \in \Z_+$, we denote by $\L_N$ the weight 2 modular form on $\G_0(N)$ given by $\L_N = N q\frac{d}{dq} (\log \frac{\h(N\t)}{\h(\t)})$ (cf. \eq{Eisenstein_form}). For $N=23$, there are two newforms and we choose the basis $f_{23,1} = \h^2_{23\!A\!B}$ and $f_{23,2}$ given in \eq{phi232}.
}}
  \end{table}

\begin{prop}\label{H_g_def}
Let $H: \H \to \C$ be given by
\be\label{second_def_H}
H(\t) = \frac{-2 E_2(\t) + 48 F_2^{(2)}(\t)}{\h(\t)^3}=2 q^{-\frac{1}{8}} \left( -1 + 45 q+ 231 q^2 +\dots \right)
\ee
where 
$$
F_2^{(2)}(\t)= \sum_{\substack{r>s>0\\ r-s=1\, {\rm mod }\; 2 }} (-1)^{r} \,s \,q^{rs/2} = q+q^2-q^3+q^4+\dots \;.
$$
Then for all $g\in M_{24}$, the function 
\be\label{h_g_explicit}
H_g(\t) = \frac{\chi(g)}{24} H(\t) - \frac{\til T_g(\t)}{\h(\t)^3}\;,
\ee
is a (mock) modular form for $\Gamma_0(N_g)$ of weight 1/2 with shadow $\chi(g)\h(\t)^3$. Moreover, 
we have
$$  
\hat H_g(\t)=\psi(\g) \, {\rm jac}(\g,\t)^{1/4}\,\hat H_g(\g\t),
$$
for $\g\in \Gamma_0(n_g)$ where 
$$
\hat{H}_g(\t)=H_g(\t)+\chi(g)\, (4{i})^{-1/2} \int_{-\bar \t}^{\infty}(z+\t)^{-1/2}\overline{\eta(-\bar z)^3}{\rm d}z.
$$
and the multiplier system is given by $\psi(\g)=\e(\g)^{\mathsmaller{-3}}\r_{n_g|h_g}(\g) $.
\end{prop}
(See \cite{Atish_Sameer_Don} for more information regarding the particular expression for $H(\t)$ given in (\ref{second_def_H}).)

Notice that the extra multiplier $\r_{n|h}$ that appears when $h_g\neq 1$ is the same as that of the inverse eta-product $1/\h_g(\t)$ (cf. \eq{rho_multiplier}) which are also related to $M_{24}$, a fact that is in accordance with the $1/2$- and $1/4$-BPS spectrum of the ${\cal N}=4$, $d=4$ theory obtained by $K3\times T^2$ compactification of the type II string theory, as will be discussed in more detail in \S\ref{StringThy}. 

Our discussion above leads to the following conjecture.
\begin{conj}\label{conj_mock}
The weight $1/2$ (mock) modular forms $H_g$ defined in \eq{h_g_explicit} satisfy 
\be\label{rep_mock}
H_g(\t) =q^{-\frac{1}{8}}\big( -2  + \sum_{n=1}^\inf q^n\,(\tr_{K_n} g)\big)
\ee
for a certain $\Z$-graded, infinite-dimensional $M_{24}$ module $K= \bigoplus_{n=1}^\inf  K_n$.

Moreover, the representations $K_n$ are even in the sense that they can all be written in the form $K_n=k_n\oplus k_n ^{\,\ast}$ for some $M_{24}$-modules $k_n$ where $k_n ^{\,\ast}$ denotes the module dual to $k_n$. 
\end{conj} 
The first few Fourier coefficients of the $q$-series $H_g(\t)$ and the corresponding $M_{24}$-representations are given in Table \ref{Mock_decompositions}. 

A proof of for the first part of the above conjecture, namely the existence of an $M_{24}$-module $K=\bigoplus_{n=1}^\inf K_n $ such that \eq{rep_mock} holds, has been attained very recently \cite{Gannon_proof}. However, it is important to stress that, unlike the eta-product moonshine described in \S\ref{Cusp}, the nature and origin of the $M_{24}$-module $K$ remains mysterious. 

The compelling modular properties of the functions $H_g$ constitute strong evidence for the existence of the $M_{24}$-module $K$. Yet stronger evidence would be furnished by a uniform construction of the (mock) modular forms $H_g$. Such a construction was established recently in \cite{Cheng2011}. The construction presented there, in terms of Rademacher sums, is of further interest in that it points to an extension of the connection between $M_{24}$ and $K3$ surfaces, expressed in CFT terms via the $H_g$ (cf. \S\ref{Weak}), to a connection in the context of quantum gravity. Also, by comparison with the Rademacher sum expressions for the functions of monstrous moonshine \cite{Duncan2009}, the Rademacher sum construction of the $H_g$ suggests a reformulation of the crucial genus zero property of monstrous moonshine that incorporates, and emphasises the distinguished nature of, the functions $H_g$ attached to $M_{24}$. 
\begin{thm} \cite{Cheng2011}
Let $g\in M_{24}$. Let $\reg(\g,\t)$ be the regularisation factor given by 
$\reg(\g,\t)=1$ in case $\g$ is upper triangular (i.e. $\g\cdot\infty=\infty$) and
\be\label{defn:reg}
\reg(\g,\t)=\ex(\tfrac{\g\t-\g\infty}{8})\,\ex(\tfrac{\g\infty-\g\t}{8},\tfrac{1}{2})
\ee
otherwise, where $\ex(x,s)$ is the following generalisation of the exponential function:
\be\label{defn:genexp}
\ex(x,s)=\sum_{m\geq 0}\frac{(\tpi x)^{m+s}}{\G(m+s+1)}\;.
\ee
Given $\G<\SL_2(\Z)$ containing the group $\G_{\infty}$ of upper triangular matrices in $\SL_2(\Z)$, and given a character $\rho$ on $\G$ define the Rademacher sum 
$$
R_{\G,\rho}(\t)
=
\lim_{K\to \inf}\sum_{\g\in(\G_{\infty}\backslash\G)_{<K}}
	\psi(\g)
	\ex(-\tfrac{\g\t}{8})
	\reg(\g,\t)
	\jac(\g,\t)^{1/4}
$$
where $\psi(\g)=\r_{n_g|h_g}(\g)\e(\g)^{\mathsmaller{-3}}$ and $(\G_{\infty}\backslash\G)_{<K}$ denotes the set of cosets of $\G_{\infty}$ in $\G$ having a representative with lower row $(c,d)$ satisfying the bounds $0\leq c<K$ and $|d|<K^2$. Then $H_g(\t)= -2 R_{\G,\rho}(\t)$ when $\G=\G_0(n_g)$ and $\rho=\r_{n_g|h_g}$. 
\end{thm}

\section{Weak Jacobi Forms and $M_{24}$} 
\label{Weak}
\setcounter{equation}{0}

In the previous section we have discussed the relation between $M_{24}$ and a set of (mock) modular forms. 
The origin and the structure of the corresponding $M_{24}$-module $K$ remains to be explored, in contrast to the situation of the cusp forms of \S\ref{Cusp}, whose relation to $M_{24}$ can be fully understood in terms of the 24-dimensional permutation representation $R$. 
But so far, the observation of the relation between $M_{24}$ and mock modular forms was in fact mostly made via a weak Jacobi form \cite{Eguchi2010}. 
This weak Jacobi form has the property that it incorporates both the infinite-dimensional $M_{24}$-module $K$ and the permutation representation $R$, a feature that will become even more significant in our later discussion in \S\ref{Siegel} about Siegel modular forms. Perhaps the best way to motivate this observation is to put it in the context of the elliptic genus of $K3$ surfaces.

For a compact complex manifold $M$ with dim$ _\C M=d_0$, the elliptic genus is defined as the character-valued Euler characteristic of the formal vector bundle 
\cite{Ochanine,Witten1987,Landweber_book,WittenInt.J.Mod.Phys.A9:4783-48001994,KawaiNucl.Phys.B414:191-2121994}
$$
{\bf E}_{q,y} 
=
y^{d/2}{\textstyle \bigwedge}{ }_{-y^{-1}} T_M^\ast
\textstyle{\bigotimes}_{n\geq 1} \textstyle\bigwedge{ }_{-y^{-1}q^n} T_M^\ast\bigotimes_{n\geq 1} \textstyle\bigwedge{ }_{-yq^n} 
T_M \bigotimes_{n\geq 0} S_{q^n} (T_M\oplus  T_M^\ast),
$$
where $T_M$ and \(T_M^\ast\) are the holomorphic tangent bundle and its dual, and we adopt the notation
\[ 
\textstyle\bigwedge{ }_q V = 1 + q V + q^2 \textstyle \bigwedge^2 V + \dots,\quad S_q V = 1 + q V +q^2 S^2V \dots,
\]
with $S^kV$ denoting the $k$-th symmetric power of $V$. In other words, we have 
\be\label{def_eg}
\EG(\t,z;M) = \int_M ch({\bf E}_{q,y}) {\rm Td}(M).
\ee
From the above definition we see that this topological quantity reduces to the familiar ones: the Euler number, the signature, and the $\hat{A}$ genus of $M$, when we specialise $z$ to $z=0,\t/2,(\t+1)/2$, respectively, when $z$ and $\t$ are such that $y=e(z)$ and $q=e(\t)$.

When $M$ has vanishing first Chern class---in particular when $M$ is a Calabi--Yau manifold---its elliptic genus $\EG(\t,z;M)$ transforms nicely under the group $\SL_2(\Z) \ltimes \Z^2$ and is in fact a weak Jacobi form of weight zero and index $d_0/2$ \cite{KawaiNucl.Phys.B414:191-2121994}. 
We say a function $\f: {\mathbb H} \times \C \to \C$ is a weak Jacobi form of weight $w$ and index $t$ if it transforms under $\SL_2(\Z) \ltimes \Z^2$ as follows:
\bea\notag
\f(\t,z) &=& \jac(\g,\t)^{\mathsmaller{w/2}} e(-t \tfrac{c z^2}{c\t+d})\, \f(\g(\t, z)) 
,\quad \g=\Big(\!\begin{array}{cc} a&b\\c&d\end{array}\!\Big) \in \SL_2(\Z)
\\\label{weak_jac}
\f(\t,z) &=&e( t(\l^2 \t + 2\l z)) \, \f(\t, z+\l \t +\m) 
,\quad \l,\,\m\in \Z\;
\eea
where 
$ \g (\t,z) = (\frac{a\t+b}{c\t+d},\frac{z}{c\t+d})$ and if the Fourier expansion of $\phi$ in $q=e(\t)$ and $y=e(z)$ involves no negative powers of $q$.

The fact that the elliptic genus of a Calabi--Yau manifold $M$ is a Jacobi form can be proven from the definition (\ref{def_eg}) \cite{Gritsenko2007}, but is especially clear from the point of view of conformal field theories. When the target space manifold $M$ is Calabi--Yau the corresponding supersymmetric sigma model with ${\cal N}=2$ supersymmetry is also conformal. The spectrum of the supersymmetric states is encoded in the function
\be\label{elliptic_genus}
\EG(\t,z;M) = \tr_{{\cal H}_{\text {RR}}}\Big( (-1)^{J_0+\bar J_0} y^{J_0} q^{L_0-c/24}\bar q^{\bar L_0-c/24}  \Big)
\ee
where $c=3d_0$ is the central charge of the conformal field theory, and $J_0, L_0,\bar J_0, \bar L_0$ are the zero modes of the left- and right-moving copies of the $\U(1)$ and Virasoro parts of the ${\cal N}=2$ superconformal algebra, respectively, and ${\cal H}_{RR}$ denotes the space of quantum states in the Ramond--Ramond sector. The action of the ${\cal N}=(2,2)$ superconformal algebra then predicts that the right hand side of (\ref{elliptic_genus}) transforms as in (\ref{weak_jac}) with $w=0$ and $t=d_0/2$. As our notation suggests, the Jacobi forms attached to $M$ by (\ref{def_eg}) and (\ref{elliptic_genus}) coincide. To understand this we should view the vector bundle ${\bf E}_{q,y}$ as the space of the Ramond--Ramond quantum states of the CFT built by acting on the Ramond sector ground states with the left-moving bosonic operators $\a_{-n}$ and the fermionic operators $\psi_{-n}^i$, $\psi_{-n}^{\bar \imath}$.  

When written in the form (\ref{elliptic_genus}) the elliptic genus may be interpreted as computing the graded dimension of the cohomology of the right-moving  supercharge  \cite{WittenInt.J.Mod.Phys.A9:4783-48001994}, graded by the left-moving quantum numbers. At the special value $z=0$ the above definition reduces to that of a Witten index enumerating Ramond--Ramond ground states. 

The Jacobi form property greatly constrains the possibilities for the elliptic genus of a Calabi--Yau manifold $M$. For example, when $M$ is a $K3$ surface the Euler number $\EG(\t,z=0;K3)=24$ forces the elliptic genus to be
\be\label{K3_EG}
	{\cal Z}(\t, z) = \EG(\t,z;K3) 
	= 8 \sum_{i=2,3,4}\left(\frac{\th_i(\t,z)}{\th_i(\t,0)}\right)^2\;,
\ee 
as computed in \cite{Eguchi1989} (cf. \S\ref{sec:JacTheta} for the definition of $\th_i(\t,z)$). In particular this function is independent of the choice of $K3$ surface and so there is no ambiguity in writing $K3$ in place of $M$ in (\ref{K3_EG}).

It is to be expected that some further special properties should hold for the elliptic genus when $M$ is a $K3$ surface because two (complex) dimensional Calabi--Yau manifolds are not only K\"ahler but also hyper-K\"ahler. As a result, the $\U(1)_R$ symmetry can be extended to $\SU(2)_R$ and the sigma model has enhanced ${\cal N}\!=4$ superconformal symmetry for both the left and the right movers. This leads to a specific decomposition of the elliptic genus of an ${\cal N}=4$ SCFT that we will now explain. 

Since the underlying CFT admits an action by the ${\cal N}=4$ superconformal algebra the Hilbert space decomposes into a direct sum of (unitary) irreducible representations of this algebra. Hence the elliptic genus (\ref{elliptic_genus}) can be written as a sum of characters of representations of this algebra with some multiplicities (cf. (\ref{K3_EG_Fou})). A natural embedding of the $\U(1)$ current algebra of the ${\cal N}=2$ superconformal algebra into the $\SU(2)$ current algebra of the ${\cal N}=4$ superconformal algebra is obatined by choosing $J_0^3 \sim J_0$. As a result, the ${\cal N}=4$ highest weight representations are again labeled by two quantum numbers $h,\ell$, corresponding to the operators $L_0$, $J_0^3$ respectively, and the character of an irreducible representation $V_{h,\ell}$ say is defined as
$$
ch_{h,\ell}(\t,z) = \tr_{V_{h,\ell}} \left( (-1)^{J_0}y^{J_0} q^{L_0-c/24}\right) \;.
$$
For the central charge $c=6$, there are two supersymmetric (also called `BPS' or `massless') representations in the Ramond sector and they have the quantum numbers 
$$
h=\tfrac{1}{4}\;,\; \ell = 0 ,\tfrac{1}{2} \;.
$$
Their characters  are given by \cite{Eguchi1987,Eguchi1988,Eguchi1988a} 
\bea\label{massless_characters}
ch_{\frac{1}{4},0}(\t,z) &=& \frac{\th_1^2(\t,z)}{\eta^3(\t)} \m(\t,z),\\\notag
ch_{\frac{1}{4},\frac{1}{2}}(\t,z) &=& q^{-\frac{1}{8}} \frac{\th_1(\t,z)}{\eta^3(\t)} - 2 \frac{\th_1(\t,z)^2}{\eta^3(\t)} \m(\t,z),
\eea 
where $\m(\t,z)$ denotes the so-called Appell-Lerch sum
 $$
 \m(\t,z) = \frac{-i y^{1/2}}{\th_{1}(\t,z)}\,\sum_{\ell=-\inf}^\inf \frac{(-1)^{\ell} y^\ell q^{\ell(\ell+1)/2}}{1-y q^\ell}.
 $$
Notice that the supersymmetric representations have non-vanishing Witten index 
$$ch_{\frac{1}{4},0}(\t,z=0) =1, \quad ch_{\frac{1}{4},0}(\t,z=0) =-2.$$ 
On the other hand the massive (or `non-BPS' or `non-supersymmetric') representations with 
$$
h =\tfrac{1}{4} + n,\quad \ell =\tfrac{1}{2},\quad\text{} \quad n = 1,2,\dots
$$
have the character given by
$$
ch_{\frac{1}{4}+n,\frac{1}{2}}(\t,z) = q^{-\frac{1}{8}+n} \frac{\th_1^2(\t,z)}{\eta^3(\t)} 
$$
which has, by definition, vanishing Witten index $ch_{\frac{1}{4}+n,\frac{1}{2}}(\t,z=0) = 0$.

Rewriting the $K3$ elliptic genus in terms of these characters, we arrive at the following specific expression for the weak Jacobi form ${\cal Z}(\t,z)$ \cite{Eguchi1989,Eguchi2008,Eguchi2009a}
\be\label{expand_mu}
{\cal Z}(\t, z) = \frac{\th_1^2(\t,z)}{\eta^3(\t)} \,\left(a \,\m(\t,z) +q^{-1/8}\big(b  + \sum_{n=1}^\inf t_n q^n \big) \right)
\ee
where $a,b\in\Z$ and $t_n\in\Z$ for all positive integers $n$, and
the $t_n$ count the number of non-supersymmetric representations with $h=n+1/4$ which contribute to the elliptic genus. 
As the notation suggests, from \eq{K3_EG} we can compute the above integers to be $a=24, b=-2$, and the first few $t_n$'s are indeed as we have seen in the last section
$$
2\times45,\, 2\times231,\, 2\times770,\,2\times2277,\,2\times 5796,\ldots 
$$
In other words, the mock modular form $H(\t)$ (cf. \eq{second_def_H}) can be interpreted as the generating function of multiplicities of massive irreducible representations of the ${\cal N}=4$ superconformal algebra in the $K3$ elliptic genus.

In this context, the origin of the mock modularity of $H(\t)$ can be understood in the following way. From the fact that the elliptic genus \eq{expand_mu} transforms nicely under $\SL_2(\Z)$, so must the combination $24 \m(\t,z)+H(\t)$. Now, the Appell-Lerch sum itself does not transform nicely, rather its non-holomorphic completion 
$$
\hat \m (\t; z) =  \m (\t; z) -\frac{1}{2} \int_{-\bar \t}^\infty dz \, \frac{{\eta^3(- z)}}{\sqrt{i(z+\tau)}}
$$
transforms like a weight $1/2$ theta function. This has been demonstrated in \cite{zwegers} as part of a systematic treatment of mock $\th$-functions. Therefore, the mock modularity of the $q$-series $H(\t)$ is directly related to the mock modularity of the massless (BPS) ${\cal N}=4$ characters \eq{massless_characters}.

In \eq{expand_mu} we have seen the appearance of the infinite-dimensional $M_{24}$-module $K$ in the elliptic genus ${\cal Z}(\t,z)$. One might wonder whether $M_{24}$ acts on the other part of the decomposition as well. A simple observation is that $a=24$ is the dimension of the defining permutation representation $R$ of $M_{24}$, and hence a naive guess will be that $M_{24}$ acts on the massless ${\cal N}=4$ multiplets as the direct sum of $R$ and a two dimensional (odd) trivial representation and this suggests that we write $24= \tr_R {\mathds 1}$. Together with this assumption, the conjecture \eq{conj_mock} of the previous section implies that $M_{24}$ acts on the states of the $K3$ sigma model contributing to the elliptic genus and moreover commutes with the superconformal algebra. Therefore, it is natural to consider also the twisted (or equivariant) elliptic genus which is expected to have the decomposition
\begin{gather}\label{K3_EG_Fou}
	\begin{split}
	 {\cal Z}_g(\t,z) 
 		&=  \frac{\th_1^2(\t,z)}{\eta^3(\t)} \,\big(\chi(g)  \,\m(\t,z) + H_g(\t)\big)\\  
 		&=\sum_{n\geq 0, \ell \in \Z} c_g(4n-\ell^2) q^n y^\ell \quad {\rm with}\;\;\chi(g)= \tr_R {g }.
	\end{split}
\end{gather}

Again, a non-trivial connection between weak Jacobi forms and $M_{24}$ arises if all such ${\cal Z}_g(\t,z)$ display interesting modular properties. 
Moreover, as we mentioned in \S \ref{Cusp}, consistency with the CFT interpretation requires this to be true. More specifically it requires that ${\cal Z}_g(\t,z)$ transform nicely under the action of $\Gamma_0(n_g)$. Indeed, from the (mock) modularity of $H_g(\t)$ (cf. \eq{H_g_def}) it is now easy to show 

\begin{prop}
For all $g\in M_{24}$, the function 
\be \label{twisted_elliptic_genus2}
{\cal Z}_g(\t,z)= \frac{\chi(g) }{12}\, \varphi_{0,1}(\t,z) + \til T_g(\t) \, \varphi_{-2,1}(\t,z)\;\ee 
is a weak Jacobi form of weight 0 and index 1 for the group $\Gamma_0(N_g)$. Moreover, we have 
$$  
{\cal Z}_g(\t,z)=\r_{n_g|h_g}(\g) e\left(-\tfrac{c z^2}{c\t+d}\right) {\cal Z}_g(\g(\t,z))
$$
for $\g \in\Gamma_0(n_g)$.
\end{prop}
In the above $\varphi_{0,1}(\t,z)$ and $\varphi_{-2,1}(\t,z)$ are weight 0 and weight 2, index 1 weak Jacobi forms whose explicit expressions can be found in Appendix \ref{Modular Forms}. The weight two modular forms $\til T_g(\t)$ are listed in Table \ref{h_g}.

As for $M_{24}$-modules, the conjecture \eq{conj_mock} on the relationship between $H_g(\t)$ and $M_{24}$ now implies that the Fourier coefficients $c_g(4n-\ell^2)$ of the twisted $K3$ elliptic genus encode the supercharacter of a $\Z$-graded super (or virtual) representation $\hat K=\bigoplus \hat K_{4n-\ell^2}$ for $M_{24}$. More specifically, we have the following conjecture.
\begin{conj} \label{conj_jac}
The Fourier coefficients $c_g(k)$ of 
$$
{\cal Z}_g(\t,z)= \sum_{n\geq 0, \ell \in \Z} c_g(4n-\ell^2) q^n y^\ell
$$
satisfy $c_g(k) = {\rm str}_{\hat K_k} g$ for $k\equiv 0,3\pmod 4$ and $k\geq -1$ where ${\hat K}_k$ is the degree $k$ component of an infinite-dimensional, $\Z$-graded super (or virtual) module $\hat K$ for $M_{24}$.
$$
\hat K =\left(\bigoplus_{\ell=0}^\inf \hat K_{4\ell-1}\right) \bigoplus\left(\bigoplus_{\ell=0}^\inf \hat K_{4\ell}\right).
$$
Moreover, when $k\neq 0$ the representation $\hat K_k$ is even in the sense that is has the form $\hat k_k \oplus  \hat k_k^{\,\ast}$ where $\hat k_k$ is a certain super representation of $M_{24}$ and $\hat k_k^{\,\ast}$ denotes its dual.
\end{conj}
The first few coefficients $c_g(k)$ and the corresponding $M_{24}$ super module $\hat K_k$ can be found in Tables \ref{elliptic_genus_decompositions1} and \ref{elliptic_genus_decompositions2}. 

We conclude this section with a related conjecture which seeks to identify a mechanism that will explain why the largest Mathieu group should be related to the modular objects described in this section.
\begin{conj} 
The elliptic cohomology of $K3$ surfaces furnishes an infinite-dimensional, $\Z$-graded $M_{24}$-supermodule.
More precisely, the cohomology of the right-moving supercharge as computed by the elliptic genus (\ref{elliptic_genus}) is an infinite-dimensional $\Z$-graded $M_{24}$-supermodule whose characters (the twisted elliptic genera) are given by ${\cal Z}_g(\t,z)$. 
 
In particular, at special points in the CFT moduli space where the Abelian subgroup $\langle g \rangle \subset M_{24}$ generated by an element $g\in M_{24}$ is a symmetry of the full Hilbert space (and not just the cohomology), we have 
 \be\label{Zg_trRR}
 {\cal Z}_g(\t,z)= \tr_{{\cal H}_{\text {RR}}}\Big( g\,(-1)^{J_0+\bar J_0} y^{J_0} q^{L_0-c/24}\bar q^{\bar L_0-c/24}  \Big)\;.
 \ee
\end{conj}
As evidence for this conjecture we observe that (\ref{Zg_trRR}) is known to be true in the case of conformal field theories corresponding to certain special points in the $K3$ moduli space where the group $\langle g \rangle$ is realised as hyper-K\"ahler automorphisms of the $K3$ surface. This is related to the result of Mukai \cite{Mukai} (see also \cite{Kondo}) that any finite group of hyper-K\"ahler automorphisms of a $K3$ surface can be embedded in $M_{23}$---the maximal subgroup of $M_{24}$ fixing a point in the defining permutation representation $R$. Mukai's result can be obtained by considering the relationship between the action of $M_{24}$ on a certain Niemeier lattice and the cohomology lattice of the $K3$ surface in question. Explicit computations along these lines can be found in \cite{Taormina2010}. An interesting stringy analogue of Mukai's result has been obtained in \cite{Gaberdiel2011}. For $g\in 2A,3A,4B,5A,6A,7A$, the twisted $K3$ elliptic genus has  been computed in this way in \cite{DavidJHEP0606:0642006,Govindarajan2008,Sen2010}. On the other hand, the elliptic cohomology can remain invariant even when the $K3$ surface is not invariant. In the generic situation where $g$ is a symmetry of the cohomology but not the full Hilbert space, one should adopt the above formula and first compute the cohomology by taking the trace on the right-moving sector before inserting the operator $g$.

\section{Siegel Modular Forms and $M_{24}$} 
\label{Siegel}
\setcounter{equation}{0}

In \S\ref{Weak} we have discussed the relation between $M_{24}$ an a set of weak Jacobi forms. It turns out that there is a set of Siegel modular forms that are closely related to these weak Jacobi forms via the so-called Borcherds, or exponential lift.  
As a result, we expect the relation described in \S\ref{Weak} to be translated into a relation between these Siegel modular forms and $M_{24}$. 
In particular, both the permutation representation $R$ and the infinite-dimensional representation $K$ play an evident role in defining these Siegel forms. 
It is the goal of this section to describe and discuss these connections.

First we will consider the relation between the weak Jacobi forms appearing in \S\ref{Weak} and Siegel modular forms. Recall that it is natural to consider the weak Jacobi forms ${\cal Z}_g(\t,z)$ as twisted elliptic genera of $K3$ surfaces as this interpretation makes the modularity as well as the connection to $M_{24}$ and the mock modular forms in \S\ref{Mock} manifest. 
Similarly, we can arrive at Siegel modular forms by considering the elliptic genera of the symmetric powers $S^NM=M^N/S_N$ of a Calabi--Yau manifold $M$. In \cite{DVV,Dijkgraaf1997}, an elegant expression for the generating function of the $\EG(\t,z;S^NM)$ for a $d_0$ (complex) dimensional Calabi--Yau manifold $M$ in terms of $\EG(\t,z;M)$ has been obtained by identifying such a generating function with the partition function of second-quantised strings on $S^1\times M$. 
This identification gives
\bea\label{SQEG}
\sum_{N=0}^\inf p^N \EG(\t,z;S^N M) 
= \exp\left( \sum_{N=0}^\inf p^N (T_N\EG)(\t,z;M) \right)
= \prod_{\substack{n > 0, \\ m\geq 0, \\ \ell \in \Z}} \Big(\frac{1}{1-p^n q^m y^\ell}\Big)^{c(2d_0 nm -\ell^2)}\;,
\eea
where 
$$\EG(\t,z;M) =\sum_{n,\ell} c(2d_0nm-\ell^2) q^n y^\ell 
$$
is the elliptic genus of $M$ itself and $T_N$ is the Hecke operator mapping a Jacobi form $\phi(\t,z)$ of weight zero and index $t$ to a form of weight zero and index $Nt$. Explicitly, we have
\be\label{hecke1}
(T_N \phi)(\t,z)=\frac{1}{N} \sum_{\substack{ad=N \\ b \text{ mod } d}} \phi\left(\frac{a\t+b}{d},az\right).
\ee
The right hand side of (\ref{hecke1}) has a natural interpretation as a sum over degree $N$ maps between elliptic curves $E\to E'$, mapping the two $1$-cycles of $E$ according to $(\begin{smallmatrix}a&b\\0&d\end{smallmatrix})$, and summing over $N$ gives the instanton contributions to the free energy of the free string theory. This leads to a path-integral derivation of the expression (\ref{SQEG}) for the second-quantised elliptic genus. As the elliptic genus of a single copy of $M$ is a weak Jacobi form, one might guess that this second-quantised elliptic genus, encoding elliptic genera of all the symmetric powers of $M$ through the parameter $p$, also has interesting automorphy. 
It was observed in \cite{DVV,Gritsenko2007} that, up to a factor which has the effect of symmetrising $p$ and $q$, the second-quantised elliptic genus is in fact a Siegel modular form when $M$ is Calabi--Yau.

To explain this recall that the {\em genus two Siegel upper half space} $\H_2$ may be realised concretely as the set of symmetric $2\times 2$ matrices with complex entries whose imaginary parts are positive definite.
\begin{gather}
	\H_2=\{\Omega\in M_2(\C)\mid \Omega=\Omega^t,\,\det\Im(\Omega)>0,\,\tr\Im(\Omega)>0\}
\end{gather}
We may introduce (complex) coordinates $\s$, $\t$ and $z$ on $\H_2$ by requiring that
\begin{gather}\label{Siegel:Omega_coords}
	\Omega=
	\begin{pmatrix}
		\s&z\\
		z&\t
	\end{pmatrix}
\end{gather}
and then $\s,\t\in\H$, where $\H$ is the upper half plane, and $z\in \C$ satisfies $\Im(z)^2<\Im(\s)\Im(\t)$. A function $\F_k(\O): \H_{2} \to \C$ is called a {\em Siegel modular form of weight $k$} for a group $\G^{(2)}\subset \Sp_2({\mathbb Q})$ if
\[ \F_k((A\O+B)(C\O+D)^{-1}) = (C\O+D)^k \F_k(\O) 
\]
for all 
\[ \bem A & B\\ C&D\eem \in \G^{(2)} \;.
\]
As alluded to above, after introducing a factor $A(\Omega)$, sometimes referred to as the Hodge anomaly, which depends only on the Hodge numbers of $M$, the second-quantised elliptic genus of an even dimensional Calabi--Yau manifold $M$ is a Siegel modular form of some weight for the full Siegel modular group $\Sp_2(\Z)$. More precisely, upon setting $p=e(\s)$ one obtains that
\[
\frac{1}{\F_k(\O)} = A(\Omega)
\sum_{N=0}^\inf p^N \EG(\t,z;S^N M)
\]
is a Siegel modular form of weight $-k=-\frac{1}{2} \chi_{d_0}'(M)$ for $\dim_{\C} M =2d_0$ where $$\chi_{n}'(M) = \sum_{m=1}^{2d_0} (-1)^{n+m} h^{n,m}(M).$$ 
See \cite{Gritsenko2007} for a complete description of the modular property.

In the case of interest to us, where $M$ is a $K3$ surface, we have 
\begin{gather}\label{BorLiftEGK3}
\begin{split}
\F_{10}(\O)	
	&=pqy\prod_{\substack{m,n,l\in\Z\\(m,n,l)>0}}(1-p^mq^ny^l)^{c(4mn-l^2)} \\
	&= pqy \prod_{(n,m,\ell)>0} \exp\left(-\sum_{k=1}^\inf \frac{c(4nm-\ell^2)}{k} (p^m q^n y^\ell)^k \right)
\end{split}
\end{gather}
for the reciprocal of the $A(\Omega)$-corrected second-quantised elliptic genus, where the condition denoted $(m,n,l)>0$ is that $m,n\geq 0$ and $l<0$ when $m=n=0$ and the values $c(k)$ are the Fourier coefficients of the $K3$ elliptic genus (\ref{K3_EG_Fou}). This function $\F_{10}(\Omega)$ is a cusp form of weight $10$ for $\Sp_2(\Z)$ and is typically taken as one of four standard generators of the full ring of Siegel modular forms  for $\Sp_2(\Z)$ \cite{igusa,igusa2}. 

Physically, the extra factor $A(\Omega)$ accounts for the further compactification on $K3\times T^2$ \cite{DavidJHEP0611:0722006,ShihJHEP0610:0872006}. As a result,  the Siegel modular form $1/\Phi_{10}$ has the interpretation of counting $1/4$-BPS states of type II superstring theory compactified on $K3\times T^2$ \cite{DVV}. See \cite{Sen2007c,Cheng2008b} for a review of this topic and for further references. There are interesting subtleties regarding the poles of the function $1/\Phi_{10}$ for which we refer also to \cite{SenJHEP0705:0392007,Cheng2007a,Dabholkar2008,Atish_Sameer_Don}. 

The elliptic genus plays an evident role in defining the automorphic form $\F_{10}(\Omega)$. Mathematically, this is a special case of a general method, known as the {\em Borcherds} or {\em exponential lift}, for transforming the data of a weak Jacobi form in to a Siegel modular form \cite{borcherds_singular,Grit_Nik,Gritsenko2008b}. In \S\ref{Weak} we have seen that, apart from ${\cal Z}(\t,z)$, we have in fact a weak Jacobi form ${\cal Z}_g(\t,z)$ for each conjugacy class $[g]$ of $M_{24}$, and a natural question to ask at this point is whether there are also Siegel modular forms $\F_g(\O)$ obtained from these by a similar lifting. It turns out that, both from the algebraic and geometric points of view, there are natural candidates for such functions on $\H_2$. 

From the interpretation of the Hecke operator as summing over degree $N$ maps between elliptic curves and the interpretation of twisting as changing the boundary condition in the path integral,  we see that we should replace the Hecke operator with the {\em equivariant} (or {\em twisted}) Hecke operator:
\be\label{hecke2}
(T_N {\cal Z}_g)(\t,z) =\frac{1}{N} \sum_{\substack{ad=N \\ b \text{ mod } d}} {\cal Z}_{g^a}\left(\frac{a\t+b}{d},az\right).
\ee
(See \cite{Carnahan2008,ganter1} for a more thorough geometric description of the equivariant Hecke operator.) From the above argument we see that the corresponding twisted object is \cite{Cheng2010_1}
\be\label{twisted_siegel}
\F_{g}(\O)=pqy \prod_{(n,m,\ell)>0} \exp\left(-\sum_{k=1}^\inf \frac{c_{g^k}(4nm-\ell^2)}{k} (p^m q^n y^\ell)^k \right).
\ee
Note that the above expression can again be rewritten in an infinite-product form. The functions $\F_g(\Omega)$ for $g\in 2A,3A,5A,7A$ were proposed previously as partition functions for $1/4$-BPS states twisted by a $\Z/p$ symmetry of the theory for $p=2,3,5,7$ in \cite{DavidJHEP0701:0162007,DabholkarJHEP0711:0772007,Sen2009,Sen2010,Govindarajan2009}. 
See also \cite{Govindarajan2011} for a closely related discussion. 

The zeros of a Siegel modular form constitute its {\em rational quadratic divisor}. In the case of $\F_{10}(\O)$ the zeros are the images under the Siegel modular group $\Sp_2(\Z)$ of 
\be\label{pole_Phi10}
\lim_{z\to 0} \,(2\p i z)^{\mathsmaller{-2}} \,\F_{10}(\O) = \h(\t)^{24}\h(\s)^{24}.
\ee
Interestingly, this limit recovers (two copies of) the cusp form $\h_{g}(\t)$ discussed in \S\ref{Cusp} for the case that $g$ is the identity element. In fact, this property holds more generally. In \S\ref{Weak} we have seen how the weak Jacobi forms ${\cal Z}_g(\t,z)$ incorporate both the permutation representation $R$, which is responsible for the connection between cusp forms $\h_g(\t)$ and $M_{24}$, and the infinite-dimensional $M_{24}$ representation $K$ underlying the connection to the mock modular forms $H_g(\t)$. The data of both these $M_{24}$-modules is thus contained in $\F_g(\O)$ via the lifting procedure described above. Through direct computation we see that taking the limit $z\to 0$ amounts to ``forgetting'' the data of the massive module $K$. And more generally, using \eq{twisting1} it can be shown \cite{Cheng2010_1} that 
\be\label{pole_Phig}
\lim_{z\to 0}\, (2\p i z)^{\mathsmaller{-2}} \,\F_{g}(\O) = \h_g(\t)\,\h_g(\s)
\ee
for all $g\in M_{24}$ with $\F_g(\Omega)$ defined as in \eq{twisted_siegel}. (See also \cite{Eguchi2011} for a slightly different approach.)

It is natural to ask about the automorphic properties of the $\F_g(\O)$. As $\F_{10}(\O)$---corresponding to the identity element of $M_{24}$---is a Siegel modular form for the full Siegel modular group, it is reasonable to guess that each $\F_g(\O)$ for $g\in M_{24}$ is a Siegel modular form for some congruence subgroup of $\Sp_2(\Z)$. Inspired by physical considerations to be discussed in the next section we are led to the following.
\begin{conj}\label{conj_Phig}
For $g\in M_{24}$ the function $\F_g(\O)$ is a Siegel modular form, with a possibly non-trivial multiplier, for a subgroup of $\Sp_2(\Z)$ containing $\G_0^{(2)}(N_g)$, where 
$$
\G_0^{(2)}(N) = \left\{ \bem A& B \\ C &D \eem \in \Sp_2(\Z) \big\lvert \, C = 0\;\, {\rm  mod   } \,\;N \right\}.
$$
\end{conj}
Note that $\G_0^{(2)}(N)$ is the genus two analogue of $\G_0(N)$ defined in \eq{congruence1}. As evidence in support of Conjecture \ref{conj_Phig} we remark that the functions $\F_g(\Omega)$ for $g\in 2A,3A,4B$ have been discussed in \cite{Gritsenko2008} and have been verified to transform in the required way under $\G_0^{(2)}(N_g)$ with $N_g=2,3,4$ respectively. The automorphic property of $\F_g(\Omega)$ has also been discussed for $g\in 5A, 7A$ in \cite{DavidJHEP0606:0642006}. 

We conlcude this section with some discussion of open questions. 
First, apart from the infinite products we have seen in (\ref{BorLiftEGK3}) and (\ref{twisted_siegel}) another way to construct a Siegel modular form of a given weight is to apply the {\em Maass} or {\em additive lift} to a Jacobi form of the same weight. See for instance \cite{Gritsenko2008b,Gritsenko2008}. It seems that the cusp form $\h_g(\t)$ can be used to perform an additive lift construction of $\F_g(\O)$ for some conjugacy classes of $M_{24}$ as well \cite{JatkarJHEP0604:0182006,Cheng2010_1,Eguchi2011}, however, so far there is no known formulation that applies to all elements of $M_{24}$.  

Second, recall that \eq{pole_Phi10} may be obtained as the {\em pinching limit} ($\e\to 0$) of the genus two bosonic string amplitude associated to gluing pairs of tori \cite{DVV,Tuite2001}. It would be interesting to have an analogous interpretation for general $\F_g(\O)$. We refer to \cite{DabholkarJHEP0712:0872007} for a related discussion relevant to the special case that $g \in 2A$. 

Last but not least, the weight five Siegel modular form $\D_5 = \F_{10}^{1/2}$ coincides with the denominator of a generalised Kac--Moody superalgebra \cite{Grit_Nik}. The relationship between $\F_{10}(\Omega)$, the Jacobi form ${\cal Z}(\t,z)$ and representations of $M_{24}$ (cf. \S\ref{Weak}) suggests the possibility of an action of $M_{24}$ on the roots of this or some closely related superalgebra. Since the Siegel modular form $\F_{10}(\Omega)$, together with its twists $\F_g(\Omega)$, unifies all the automorphic objects discussed above it would be very interesting to have a concrete $M_{24}$-module that realises these functions $\F_g(\Omega)$. The existence of the generalised Kac--Moody superalgebra attached to $\D_5$ suggests that such a construction might be possible. See \cite{Govindarajan2010a} for some related discussions. We mention here that early results relating Siegel modular forms to infinite-dimensional Kac--Moody algebras appeared in \cite{feingold_frenkel}. In particular, evidence for a family of generalisations of the Maass lift (mentioned above) is provided in loc. cit. According to \cite{feingold_nicolai} the subalgebra of the generalised Kac--Moody superalgebra attached to $\D_5$ generated by its real roots appears as a subalgebra of the hyperbolic Kac--Moody algebra considered in \cite{feingold_frenkel}. Interesting connections between $\D_5$, moduli of $K3$ surfaces and (an ``automorphic correction'' of) the hyperbolic Kac--Moody algebra of \cite{feingold_frenkel} are established in \cite{GriNik_IguLKM}.

\section{String Theory and $M_{24}$}\label{StringThy}

\begin{figure}
\begin{center}
\includegraphics[scale=0.4]{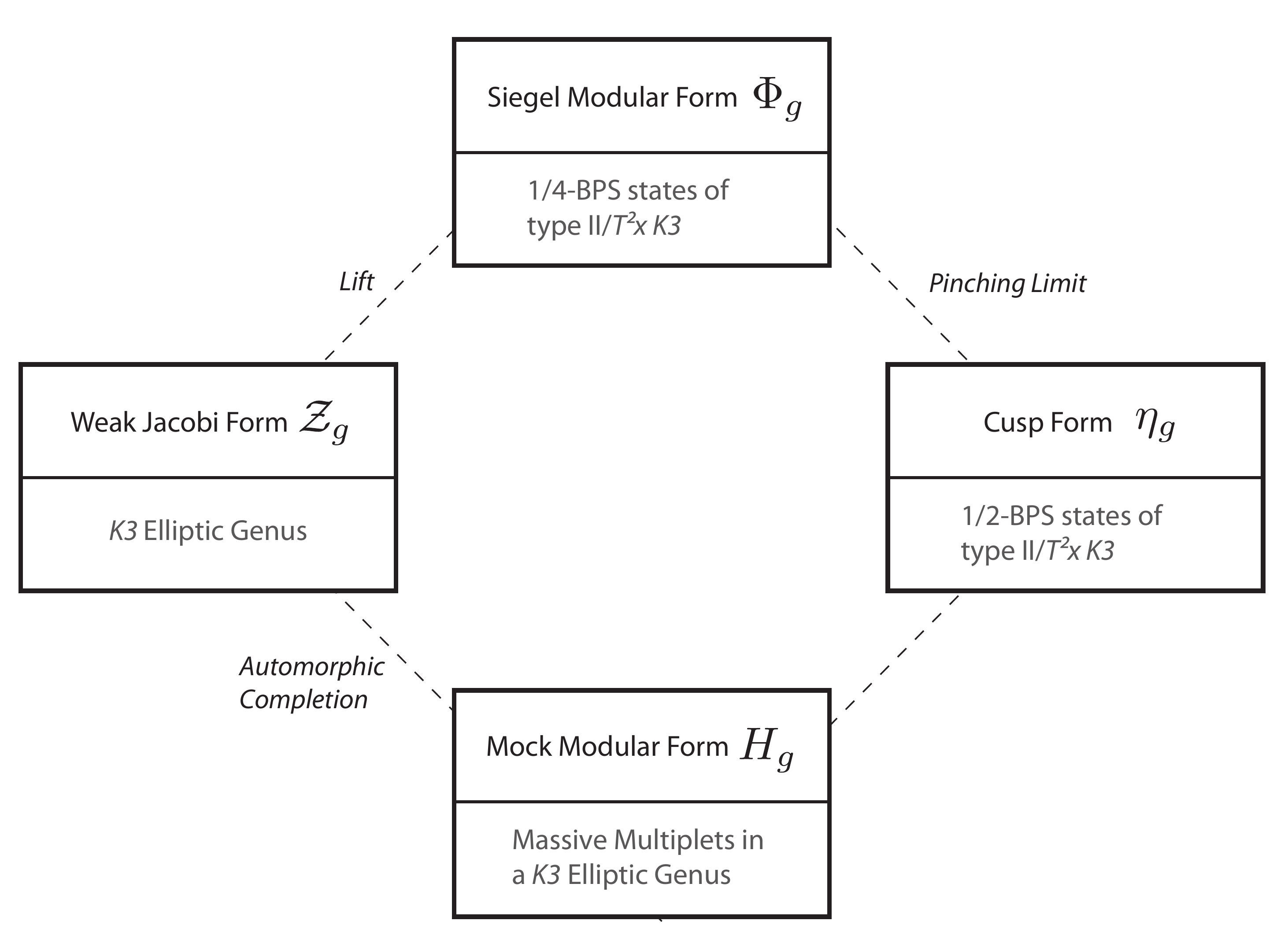}
\caption{The relations between the automorphic objects discussed in this review and the physical objects associated to them.\label{fig:diag}}
\label{diagram}
\end{center}
\end{figure}

So far we have phrased our discussions in terms of automorphic forms and the relatively familiar language of conformal field theories and vertex operator algebras. 
However, much useful intuition can be been derived from the physical system of quantum black holes in the $K3\times T^2$-compactified string theory---another indication that the largest Mathieu group is closely connected to $K3$ surfaces in a way that is yet to be fully understood.  In this final section we will describe how the various mathematical objects and relations discussed above can be understood in the context of this physical setup. We will assume familiarity with the basics of string theory. 

First we start with the string theoretic interpretation of the cusp forms of \S \ref{Cusp}. 
Apart from being the partition function of the theory with $24$ free chiral bosons, the function $1/\h(\t)^{24}$, which is $1/\h_g(\t)$ for $g$ the identity of $M_{24}$, encodes the graded numbers of the  supersymmetric (perturbative) states in heterotic string theory. 
As a result, the duality between type II and heterotic string theories dictates that this function is also the generating function for the index enumerating BPS states preserving $8$ of the 16 supercharges of type II string theory compactified on $K3\times T^2$. In the event that a particular $K3$ surface admits a hyper-K\"ahler automorphism that can be identified with the cyclic subgroup $\langle g \rangle$ of $M_{24}$ \cite{Mukai,Kondo} we can also compute the $1/2$-BPS partition function twisted by the symmetry $g$ of the Hilbert space. 
As can be expected from the discussion in \S\ref{Cusp} these turn out to be none other than inverses $Z_g(\t)= 1/\h_g(\t)$ of cusp forms attached to $M_{24}$ in the cases that have been explicitly verified; viz. $g \in 2A,3A,4B,5A$ \cite{Sen2007c,Govindarajan2009}. More generally, if the $1/2$-BPS spectrum of the $K3\times T^2$-compactified type II string theory has a symmetry generated by $g\in M_{24}$ we can expect the corresponding twisted partition function to be given by $Z_g(\t)= 1/\h_g(\t)$. 

Next we turn to the physical interpretation of the connection between representations of $M_{24}$ and the Siegel modular form $\F_{10}(\O)$ arising from the fact that the latter can be identified (cf. (\ref{BorLiftEGK3}))with the Borcherds lift of ${\cal Z}(\t,z)$. 
Recall that $\F_{10}(\O)$ is almost equal to the generating function of the elliptic genera of the symmetric powers $S^NK3$ and only differs from this by the Hodge anomaly $A(\Omega)$. Physically this extra factor $A(\Omega)$ accounts for the further compactification on $K3\times T^2$ \cite{DavidJHEP0611:0722006,ShihJHEP0610:0872006}. As a result,  the function form $1/\Phi_{10}$ can be interpreted as counting $1/4$-BPS states of type II superstring theory compactified on $K3\times T^2$ \cite{DVV}. 
(See \cite{Sen2007c,Cheng2008b} for a review of this topic and for further references.)
With this interpretation, Conjecture \ref{conj_jac} 
suggests that these $1/4$-BPS states of  $K3\times T^2$-compactified superstring theory furnish representations of $M_{24}$, which then suggests that we can compute the twisted $1/4$-BPS partition function for the theory and arrive at $\F_g(\O)$ as given in \eq{twisted_siegel}.  

It is interesting to observe that  \eq{pole_Phig}  must hold in order for the above string theoretic interpretation of $1/\Phi_g(\O)$ and $1/\h_g(\t)$ to be consistent \cite{Cheng2010_1}.  This is because the number of states preserving $1/2$ and $1/4$ of the supersymmetries are related through the presence of the so-called two-centred bound states, which are bound states preserving $1/4$ of the supersymmetries of two particles that preserve $1/2$ of the supersymmetries in their respective isolation. 
As alluded to earlier, the zeros of the Siegel modular form $\F_{10}(\O)$ play a special role in the physical counting of BPS states. 
They correspond to the poles of the partition function $1/\F_{10}$ and their physical significance is different depending on whether or not they intersect the (two-sheeted) hyperboloid $\langle\Im(\O),\Im(\O)\rangle = -R$ when $R$ is taken to be very large (and real). On the other hand, the  poles which do intersect the hyperboloid, such as those at $z=0$,  control the absence or presence of the aforementioned two-centred bound states \cite{SenJHEP0705:0392007,Cheng2007a,Dabholkar2008,Atish_Sameer_Don}. Recall that such $1/2$-BPS states can be thought of as quantum states in a theory with 24 free bosons and hence admit an action by $M_{24}$ according to the discussion of \S\ref{Cusp}; this dictates \eq{pole_Phig}. We summarise the relations between the cusp forms, weak Jacobi forms, Siegel modular forms and their physical interpretations in Figure \ref{fig:diag}.
 
The above discussion constitutes strong evidence that the largest Mathieu group acts as symmetries of the supersymmetric states of $K3\times T^2$-compactified superstring theory at any point in the moduli space of the theory. An understanding of this is yet to be achieved but some related results have been obtained recently in \cite{Taormina2010,Gaberdiel2011}. This problem is extremely interesting since a positive solution will provide us with a way to understand the relation between $M_{24}$ and the automorphic objects discussed in the present paper through string theory. 

\appendix

\section{Modular Forms}
\label{Modular Forms}

\subsection{Dedekind Eta Function}
The {\em Dedekind eta function}, denoted $\eta(\t)$, is a holomorphic function on the upper half-plane defined by the infinite product 
$$\eta(\t)=q^{1/24}\prod_{n\geq 1}(1-q^n)$$
where $q=\ex(\t)=e^{\tpi \t}$. It is a modular form of weight $1/2$ for the modular group $\SL_2(\Z)$ with multiplier $\e:\SL_2(\Z)\to\C^*$, which means that 
$$\e(\g)\eta(\g\t)\jac(\g,\t)^{1/4}=\eta(\t)$$
for all $\g = \big(\begin{smallmatrix} a&b\\ c&d \end{smallmatrix}\big) \in\SL_2(\Z)$, where $\jac(\g,\t)=(c\t+d)^{-2}$. The {\em multiplier system} $\e$ may be described explicitly as 
\be\label{Dedmult}
\e\bem a&b\\ c&d\eem 
=
\begin{cases}
	\ex(-b/24),&c=0,\,d=1\\
	\ex(-(a+d)/24c+s(d,c)/2+1/8),&c>0
\end{cases}
\ee
where $s(d,c)=\sum_{m=1}^{c-1}(d/c)((md/c))$ and $((x))$ is $0$ for $x\in\Z$ and $x-\lfloor x\rfloor-1/2$ otherwise. We can deduce the values $\e(a,b,c,d)$ for $c<0$, or for $c=0$ and $d=-1$, by observing that $\e(-\g)=\e(\g)\ex(1/4)$ for $\g\in\SL_2(\Z)$.

Let $T$ denote the element of $\SL_2(\Z)$ such that $\tr(T)=2$ and $T\t=\t+1$ for $\t\in \H$. Observe that
$$
\e(T^m\g)=\e(\g T^m)=\ex(-m/24)\e(\g)
$$
for $m\in\Z$.

\subsection{Theta Functions}\label{sec:JacTheta}
The {\em Jacboi theta functions} $\th_i(\t,z)$ are defined as follows for $q=e(\t)$ and $y=e(z)$.
$$	\th_1(\t,z)
	= -i q^{1/8} y^{1/2} \prod_{n=1}^\inf (1-q^n) (1-y q^n) (1-y^{-1} q^{n-1})$$
$$	\th_2(\t,z)
	=  q^{1/8} y^{1/2} \prod_{n=1}^\inf (1-q^n) (1+y q^n) (1+y^{-1} q^{n-1})$$ 
$$	\th_3(\t,z)
	=  \prod_{n=1}^\inf (1-q^n) (1+y \,q^{n-1/2}) (1+y^{-1} q^{n-1/2})$$
$$	\th_4(\t,z) 
	=  \prod_{n=1}^\inf (1-q^n) (1-y \,q^{n-1/2}) (1-y^{-1} q^{n-1/2})$$

\subsection{Weak Jacobi Forms}

According to \cite[Thm. 9.3]{eichler_zagier} (see also \cite[\S7]{feingold_frenkel}) the ring of weak Jacobi forms of even weight is generated by the functions
\bea
\notag
\varphi_{0,1}(\t,z) &=& 4\,\left(\left(\frac{\th_2(\t,z)}{\th_2(\t,0)}\right)^2 +\left(\frac{\th_3(\t,z)}{\th_3(\t,0)}\right)^2 +\left(\frac{\th_4(\t,z)}{\th_4(\t,0)}\right)^2 \right),\\\notag
\varphi_{-2,1}(\t,z)&=& -\frac{\th_1(\t,z)^2}{\eta(\t)^6},
\eea
where $\varphi_{0,1}$ has weight $0$ and index $1$ and $\varphi_{-2,1}$ has weight $-2$ and index $1$.

\subsection{Forms of Higher Level}

The congruence subgroups of the modular group $\SL_2(\Z)$ that are most relevant for this paper are 
\bea\notag
\G_0(N) &=& \bigg\{\bigg[\begin{array}{cc} a&b\\c&d\end{array}\bigg]\in \SL_2(\Z)  , c= 0 \text{ mod }N\;
\bigg\}.
\eea
For $N>1$ a (non-zero) modular form of weight two for $\G_0(N)$ is given by
\bea\label{Eisenstein_form}
\L_N(\t)&=&N\, q\pa_q\log\left(\frac{\eta(N\tau)}{\eta(\tau)}\right)\\\notag&=&\frac{N(N-1)}{24}\left(1+\frac{24}{N-1}\sum_{k>0}\s(k) (q^k -N q^{Nk})\right)
\eea
where $\s(k)$ is the divisor function $\s(k)=\sum_{d\lvert k}d$.

A modular form on $\Gamma_0(N)$ is a modular form on $\Gamma_0(M)$ whenever $N|M$ and for some small $N$ the space of forms of weight $2$ spanned by the $\L_d(\t)$ for $d$ a divisor of $N$. In the case that $N=11$ we have the {\em newform} 
\be\label{11newforms}
f_{11}(\t) = \eta^{2}(\t)\eta^{2}(11\t)
\ee
which is cusp form of weight $2$ for $\Gamma_0(11)$ that is not a multiple of $\L_{11}(\t)$. We meet the newforms
\bea
f_{14}(\t) &=& \eta(\t)\eta(2\t)\eta(7\t)\eta(14\t),\\
f_{15}(\t) &=& \eta(\t)\eta(3\t)\eta(5\t)\eta(15\t),
\eea
at $N=14$ and $N=15$, respectively, and together with $f_N$ the functions $\L_d(\t)$ for $d|N$ span the space of weight $2$ forms on $\G_0(N)$ for $N=11,14,15$.

For $N=23$ there is a two dimensional space of newforms. We may use the basis
\bea\notag
 f_{23,1}(\t)&=&\h_{23AB}(\t)^2= \h(\t)^2\h(23\t)^2\\ \label{phi232}
f_{23,2}(\t)&=& \frac{\eta(\tau)^3\eta(23\tau)^3}{\eta(2\tau)\eta(46\tau)}
	+4\eta(\tau)\eta(2\tau)\eta(23\tau)\eta(46\tau)
	+4\eta(2\tau)^2\eta(46\tau)^2
\eea
(but note that these are not Hecke eigenforms). See \cite{modi} and Chapter 4.D of \cite{from_number} for more details.

A discussion of the ring of weak Jacobi forms of higher level can be found in \cite{aoki}. From its Proposition 6.1 we conclude that the space of weak Jacobi forms of weight $0$ and index $1$ for $\G_0(N)$ is spanned by $\varphi_{0,1}(\t,z)$ and the functions $f(\t)\varphi_{-2,1}(\t,z)$ where $f(\t)$ is a modular form of weight $2$ for $\G_0(N)$.
\vspace{2cm}

 \begin{sidewaystable} 
 \section{Character Table} 
\label{Character Tables}\centering
 \resizebox{1.0\textwidth}{!}{
 \begin{tabular}[H]{ccccccccccccccccccccccccccc}
 \toprule
classes&$1A$&$2A$&$2B$&$3A$&$3B$&$4A$&$4B$&$4C$&$5A$&$6A$&$6B$&$7A$&$\overline{7A}$&$8A$&$10A$&$11A$&$12A$&$12B$&$14A$&$\overline{14A}$&$15A$&$\overline{15A}$&$21A$&$\overline{21A}$&$23A$&$\overline{23A}$\\
\midrule
 &1 & 1 & 1 & 1 & 1 & 1 & 1 & 1 & 1 & 1 & 1 & 1 & 1 & 1 & 1 & 1 & 1 & 1 & 1 & 1 & 1 & 1 & 1 & 1 & 1 & 1\\
&23&7&-1&5&-1&-1&3&-1&3&1&-1&2&2&1&-1&1&-1&-1&0&0&0&0&-1&-1&0&0\\&
45 & - 3 & 5 & 0 & 3 & - 3 & 1 & 1 & 0 & 0 & - 1 & $e_7$ & $\bar e_7$& - 
                    1 & 0 & 1 & 0 & 1 & - $e_7$ & - 
           $\bar e_7$ & 0 & 0 & $e_7$ & $\bar e_7$& - 1 & - 1\\&
           ${\overline{45}}$&-3&5&0&3&-3&1&1&0&0&-1&$\bar e_ 7$&$e_ 7$&-1&0&1&0&1&-$\bar e_ 7$&-$e_ 7$&0&0&$\bar e_ 7$&$e_ 7$&-1&-1\\&
                   231& 7& -9& -3& 0& -1& -1& 3& 1& 1& 0& 0& 0& -1& 1& 0& -1& 0& 0& 0& 
$e_{15}$& $\bar e_{15}$& 0& 0& 1& 1\\&
                  $\overline{231}$&7&-9&-3&0&-1&-1&3&1&1&0&0&0&-1&1&0&-1&0&0&0&$\bar e_{15}$&$e_{15}$&0&0&1&1\\&
                    252&28&12&9&0&4&4&0&2&1&0&0&0&0&2&-1&1&0&0&0&-1&-1&0&0&-1&-1\\&
                 253&13&-11&10&1&-3&1&1&3&-2&1&1&1&-1&-1&0&0&1&-1&-1&0&0&1&1&0&0\\&
                    483&35&3&6&0&3&3&3&-2&2&0&0&0&-1&-2&-1&0&0&0&0&1&1&0&0&0&0\\&
                    770&-14&10&5&-7&2&-2&-2&0&1&1&0&0&0&0&0&-1&1&0&0&0&0&0&0&$e_{23}$&$\bar e_{23}$\\&
                    $\overline{770}$&-14&10&5&-7&2&-2&-2&0&1&1&0&0&0&0&0&-1&1&0&0&0&0&0&0&$\bar e_{23}$&$e_{23}$\\&
                    990&-18&-10&0&3&6&2&-2&0&0&-1&$e_ 7$&$\bar e_ 7$&0&0&0&0&1&$e_ 7$&$\bar e_ 7$&0&0&$e_ 7$&$\bar e_ 7$&1&1\\& 
                    $\overline{990}$&-18&-10&0&3&6&2&-2&0&0&-1&$\bar e_ 7$&$e_ 7$&0&0&0&0&1&$\bar e_ 7$&$e_ 7$&0&0&$\bar e_ 7$&$e_ 7$&1&1\\&
                    1035&27&35&0&6&3&-1&3&0&0&2&-1&-1&1&0&1&0&0&-1&-1&0&0&-1&-1&0&0\\&
                    1035&-21&-5&0&-3&3&3&-1&0&0&1&2$e_ 7$&2$\bar e_ 7$&-1&0&1&0&-1&0&0&0&0&-$e_ 7$&-$\bar e_ 7$&0&0\\&
                    $\overline{1035}$&-21&-5&0&-3&3&3&-1&0&0&1&2$\bar e_ 7$&2$e_ 7$&-1&0&1&0&-1&0&0&0&0&-$\bar e_ 7$&-$e_ 7$&0&0\\&
                  1265&49&-15&5&8&-7&1&-3&0&1&0&-2&-2&1&0&0&-1&0&0&0&0&0&1&1&0&0\\& 1771&-21&11&16&7&3&-5&-1&1&0&-1&0&0&-1&1&0&0&-1&0&0&1&1&0&0&0&0\\& 2024&8&24&-1&8&8&0&0&-1&-1&0&1&1&0&-1&0&-1&0&1&1&-1&-1&1&1&0&0\\& 2277&21&-19&0&6&-3&1&-3&-3&0&2&2&2&-1&1&0&0&0&0&0&0&0&-1&-1&0&0\\& 3312&48&16&0&-6&0&0&0&-3&0&-2&1&1&0&1&1&0&0&-1&-1&0&0&1&1&0&0\\& 3520&64&0&10&-8&0&0&0&0&-2&0&-1&-1&0&0&0&0&0&1&1&0&0&-1&-1&1&1\\& 5313&49&9&-15&0&1&-3&-3&3&1&0&0&0&-1&-1&0&1&0&0&0&0&0&0&0&0&0\\& 5544&-56&24&9&0&-8&0&0&-1&1&0&0&0&0&-1&0&1&0&0&0&-1&-1&0&0&1&1\\& 5796&-28&36&-9&0&-4&4&0&1&-1&0&0&0&0&1&-1&-1&0&0&0&1&1&0&0&0&0\\& 10395&-21&-45&0&0&3&-1&3&0&0&0&0&0&1&0&0&0&0&0&0&0&0&0&0&-1&-1\\
                   \bottomrule
  \end{tabular}}
  \caption{ \label{M24}\footnotesize{Character table of $M_{24}$. See \cite{atlas}. We adopt the naming system of \cite{atlas} and use the notation $e_n = \frac{1}{2}(-1+i\sqrt{n})$. }}
 \end{sidewaystable}

 \begin{sidewaystable} \section{Tables of Decompositions}
\label{Tables of Fourier Coefficients and Decompositions of Representations}
\centering
 \resizebox{1.08\textwidth}{!}{
 \begin{tabular}[H]{ccccccccccccccccccccccccccc}
 \toprule
$1/\eta_g(\t)$&$1A$&$2A$&$2B$&$3A$&$3B$&$4A$&$4B$&$4C$&$5A$&$6A$&$6B$&$7A$&$\overline{7A}$&$8A$&$10A$&$11A$&$12A$&$12B$&$14A$&$\overline{14A}$&$15A$&$\overline{15A}$&$21A$&$\overline{21A}$&$23A$&$\overline{23A}$\\\midrule
$q^{0}$&24&8&0&6&0&0&4&0&4&2&0&3&3&2&0&2&0&0&1&1&1&1&0&0&1&1\\
$q^{1}$&324&52&12&27&0&4&16&0&14&7&0&9&9&6&2&5&1&0&3&3&2&2&0&0&2&2\\
$q^{2}$&3200&256&0&104&8&0&48&0&40&16&0&22&22&12&0&10&0&0&4&4&4&4&1&1&3&3\\
$q^{3}$&25650&1122&90&351&0&18&142&6&105&39&0&51&51&28&5&20&3&0&9&9&6&6&0&0&5&5\\
$q^{4}$&176256&4352&0&1080&0&0&368&0&256&80&0&108&108&52&0&36&0&0&12&12&10&10&0&0&7&7\\
$q^{5}$&1073720&15640&520&3107&44&56&928&0&590&175&4&221&221&104&10&65&5&0&23&23&17&17&2&2&11&11\\
$q^{6}$&5930496&52224&0&8424&0&0&2176&0&1296&336&0&432&432&184&0&110&0&0&32&32&24&24&0&0&15&15
\\$q^{7}$&30178575&165087&2535&21762&0&175&4979&27&2740&666&0&819&819&341&20&185&10&0&55&55&37&37&0&0&22&22
\\$q^{8}$&143184000&495872&0&53976&192&0&10864&0&5600&1232&0&1506&1506&580&0&300&0&0&76&76&56&56&3&3&30&30
\\$q^{9}$&639249300&1428612&10908&129141&0&468&23184&0&11130&2289&0&2706&2706&1010&38&481&15&0&122&122&81&81&0&0&42&42\\
                \bottomrule
                &&&& &&&& &&&& &&&&&&&& &&&&&&\\
                &&&& &&&& &&&& &&&&&&&& &&&&&&\\
                &&&& &&&& &&&& &&&&&&&& &&&&&&\\
                \toprule
                &$ \rho_{1}$&$ \rho_{23}$&$ \rho_{45}$&$ \rho_{\overline{45}}$&$ \rho_{231}$&$ \rho_{\overline{231}}$&$ \rho_{252}$&$ \rho_{253}$&$ \rho_{483}$&$ \rho_{770}$&$ \rho_{\overline{770}}$&$ \rho_{990}$&$ \rho_{\overline{990}}$&$ \rho_{1035'}$&$ \rho_{1035}$&$ \rho_{\overline{1035}}$&$ \rho_{1265}$&$ \rho_{1771}$&$ \rho_{2024}$&$ \rho_{2277}$&$ \rho_{3312}$&$ \rho_{3520}$&$ \rho_{5313}$&$ \rho_{5544}$&$ \rho_{5796}$&$ \rho_{10395}$ \\\midrule
               ${\cal H}_1$ &1&1&0&0&0&0&0&0&0&0&0&0&0&0&0&0&0&0&0&0&0&0&0&0&0&0\\
               ${\cal H}_2$ &3&3&0&0&0&0&1&0&0&0&0&0&0&0&0&0&0&0&0&0&0&0&0&0&0&0\\
               ${\cal H}_3$ &6&8&0&0&0&0&3&2&1&0&0&0&0&0&0&0&1&0&0&0&0&0&0&0&0&0\\
              ${\cal H}_4$  &14&20&0&0&0&0&12&6&4&0&0&0&0&1&0&0&3&0&0&0&1&3&0&0&0&0\\
               ${\cal H}_5$ &27&48&0&0&0&0&33&22&15&0&0&0&0&3&0&0&15&1&0&3&6&16&3&0&0&3\\
                ${\cal H}_6$&59&110&0&0&0&0&97&61&51&0&0&0&0&19&0&0&54&10&9&17&34&69&25&6&9&26\\
               ${\cal H}_7$& 114&249&0&0&6&6&255&174&161&3&3&3&3&70&3&3&190&45&47&88&158&276&147&59&71&194
               \\${\cal H}_8$&235&552&0&0&32&32&687&457&498&40&40&39&39&301&39&39&633&220&269&393&694&1042&758&418&490&1088\\${\cal H}_9$&460&1217&1&1&169&169&1783&1235&1504&255&255&296&296&1126&294&294&2152&994&1252&1730&2850&3870&3528&2354&2656&5544\\${\cal H}_{10}$&924&2677&24&24&731&731&4754&3294&4575&1425&1425&1675&1675&4329&1699&1699&7207&4391&5592&7131&11460&14340&15393&11758&13026&25565\\
                \bottomrule
  \end{tabular}}
  \caption{ \label{eta_prod_decomp}\footnotesize{i) The first ten Fourier coefficients of the McKay-Thompson series $Z_g(\t) = \sum_{k=0}^\inf q^{k-1} (\tr_{{\cal H}_k} g) =1/\eta_g(\t)$. ii) The decomposition of the first ten $M_{24}$-modules ${\cal H}_k$ into irreducible representations.  } }
 \end{sidewaystable}

 \begin{sidewaystable} \centering
    \resizebox{1.08\textwidth}{!}{
 \begin{tabular}[h!]{c|cccccccccccccccccccccccccc}
 \toprule
{\backslashbox{$\ell$}{$[g]$}}&$1A$&$2A$&$2B$&$3A$&$3B$&$4A$&$4B$&$4C$&$5A$&$6A$&$6B$&$7A$&$\overline{7A}$&$8A$&$10A$&$11A$&$12A$&$12B$&$14A$&$\overline{14A}$&$15A$&$\overline{15A}$&$21A$&$\overline{21A}$&$23A$&$\overline{23A}$\\\midrule
  0&2&2&2&2&2&2&2&2&2&2&2&2&2&2&2&2&2&2&2&2&2&2&2&2&2&2\\1&-128&0&0&-2&4&16&0&8&2&6&12&5&5&8&10&4&10&8&7&7&8&8&11&11&10&10\\2&-1026&-2&-2&0&0&62&-2&14&4&16&40&17&17&30&28&8&32&20&19&19&25&25&35&35&32&32\\3&-5504&0&0&4&4&176&0&24&6&36&108&47&47&88&70&18&92&48&49&49&69&69&95&95&85&85\\4&-23550&2&2&-6&12&450&2&50&10&74&260&110&110&226&162&34&234&104&114&114&169&169&236&236&209&209\\5&-86400&0&0&0&0&1072&0&88&20&144&576&239&239&536&340&60&544&208&245&245&380&380&539&539&471&471\\
    \bottomrule \multicolumn{27}{c}{}\\\multicolumn{27}{c}{}\\\multicolumn{27}{c}{} \\\toprule
    {\backslashbox{$\ell$}{$\rho$}}&$ \rho_{1}$&$ \rho_{23}$&$ \rho_{45}$&$ \rho_{\overline{45}}$&$ \rho_{231}$&$ \rho_{\overline{231}}$&$ \rho_{252}$&$ \rho_{253}$&$ \rho_{483}$&$ \rho_{770}$&$ \rho_{\overline{770}}$&$ \rho_{990}$&$ \rho_{\overline{990}}$&$ \rho_{1035'}$&$ \rho_{1035}$&$ \rho_{\overline{1035}}$&$ \rho_{1265}$&$ \rho_{1771}$&$ \rho_{2024}$&$ \rho_{2277}$&$ \rho_{3312}$&$ \rho_{3520}$&$ \rho_{5313}$&$ \rho_{5544}$&$ \rho_{5796}$&$ \rho_{10395}$ \\\midrule
 $0$&   2&0&0&0&0&0&0&0&0&0&0&0&0&0&0&0&0&0&0&0&0&0&0&0&0&0\\ 
 $1$&8&-2&-1&-1&0&0&0&0&0&0&0&0&0&0&0&0&0&0&0&0&0&0&0&0&0&0\\ 
 $2$&24&-6&-5&-5&-1&-1&0&0&0&0&0&0&0&0&0&0&0&0&0&0&0&0&0&0&0&0\\ 
 $3$&64&-16&-15&-15&-5&-5&0&0&0&-1&-1&0&0&0&0&0&0&0&0&0&0&0&0&0&0&0\\ 
 $4$&154&-40&-40&-40&-15&-15&0&0&0&-5&-5&0&0&0&0&0&0&0&0&-2&0&0&0&0&0&0\\ 
 $5$&342&-90&-97&-97&-40&-40&0&0&0&-15&-15&0&0&0&0&0&0&0&0&-10&0&0&0&0&-2&0\\ \bottomrule
     \end{tabular}}
     \caption{\label{elliptic_genus_decompositions1} \footnotesize{The first few Fourier coefficients $c_g(4\ell-1)$ and the corresponding super representation $\hat K_{4\ell-1}$  of $M_{24}$.  }}

   \vspace{2\baselineskip}
   \resizebox{1.08\textwidth}{!}{
 \begin{tabular}[h!]{c|cccccccccccccccccccccccccc}
 \toprule
{\backslashbox{$\ell$}{$[g]$}}&$1A$&$2A$&$2B$&$3A$&$3B$&$4A$&$4B$&$4C$&$5A$&$6A$&$6B$&$7A$&$\overline{7A}$&$8A$&$10A$&$11A$&$12A$&$12B$&$14A$&$\overline{14A}$&$15A$&$\overline{15A}$&$21A$&$\overline{21A}$&$23A$&$\overline{23A}$\\\midrule
0&20&4&-4&2&-4&-4&0&-4&0&-2&-4&-1&-1&-2&-4&-2&-4&-4&-3&-3&-3&-3&-4&-4&-3&-3\\ 
$1$&216&-8&8&0&0&-24&0&-8&-4&-8&-16&-8&-8&-12&-12&-4&-12&-8&-8&-8&-10&-10&-14&-14&-14&-14\\ $2$&1616&16&-16&-4&-4&-80&0&-16&-4&-20&-52&-22&-22&-40&-36&-12&-44&-28&-26&-26&-34&-34&-46&-46&-40&-40\\ $3$&8032&-32&32&4&-8&-224&0&-32&-8&-44&-136&-60&-60&-112&-88&-20&-116&-56&-60&-60&-86&-86&-120&-120&-110&-110\\ $4$&33048&56&-56&0&0&-568&0&-56&-12&-88&-320&-132&-132&-284&-196&-40&-292&-128&-140&-140&-210&-210&-294&-294&-256&-256\\ $5$&117280&-96&96&-8&-8&-1312&0&-96&-20&-168&-696&-292&-292&-656&-404&-68&-664&-240&-292&-292&-458&-458&-652&-652&-572&-572\\
    \bottomrule \bottomrule \multicolumn{27}{c}{}\\\multicolumn{27}{c}{}\\\multicolumn{27}{c}{} \\\toprule
   {\backslashbox{$\ell$}{$\r$}} &$ \rho_{1}$&$ \rho_{23}$&$ \rho_{45}$&$ \rho_{\overline{45}}$&$ \rho_{231}$&$ \rho_{\overline{231}}$&$ \rho_{252}$&$ \rho_{253}$&$ \rho_{483}$&$ \rho_{770}$&$ \rho_{\overline{770}}$&$ \rho_{990}$&$ \rho_{\overline{990}}$&$ \rho_{1035'}$&$ \rho_{1035}$&$ \rho_{\overline{1035}}$&$ \rho_{1265}$&$ \rho_{1771}$&$ \rho_{2024}$&$ \rho_{2277}$&$ \rho_{3312}$&$ \rho_{3520}$&$ \rho_{5313}$&$ \rho_{5544}$&$ \rho_{5796}$&$ \rho_{10395}$ \\\midrule
      ${0}$  &-3&1&0&0&0&0&0&0&0&0&0&0&0&0&0&0&0&0&0&0&0&0&0&0&0&0\\ 
    $1$&-10&2&2&2&0&0&0&0&0&0&0&0&0&0&0&0&0&0&0&0&0&0&0&0&0&0\\ 
    $2$&-32&8&6&6&2&2&0&0&0&0&0&0&0&0&0&0&0&0&0&0&0&0&0&0&0&0\\ 
    $3$&-80&20&20&20&6&6&0&0&0&2&2&0&0&0&0&0&0&0&0&0&0&0&0&0&0&0\\ 
    $4$&-190&50&50&50&20&20&0&0&0&6&6&0&0&0&0&0&0&0&0&4&0&0&0&0&0&0\\ 
    $5$&-412&108&120&120&50&50&0&0&0&20&20&0&0&0&0&0&0&0&0&12&0&0&0&0&4&0\\\bottomrule
     \end{tabular}}
     \caption{\label{elliptic_genus_decompositions2} \footnotesize{The first few Fourier coefficients $c_g(4\ell)$ and the corresponding super representation $\hat K_{4\ell}$  of $M_{24}$.  }}
       \end{sidewaystable}

 \begin{sidewaystable} \centering
    \resizebox{1.08\textwidth}{!}{
 \begin{tabular}[h!]{c|cccccccccccccccccccccccccc}
 \toprule
{\backslashbox{$\ell$}{$[g]$}}&$1A$&$2A$&$2B$&$3A$&$3B$&$4A$&$4B$&$4C$&$5A$&$6A$&$6B$&$7A$&$\overline{7A}$&$8A$&$10A$&$11A$&$12A$&$12B$&$14A$&$\overline{14A}$&$15A$&$\overline{15A}$&$21A$&$\overline{21A}$&$23A$&$\overline{23A}$\\\midrule
1&90&-6&10&0&6&-6&2&2&0&0&-2&-1&-1&-2&0&2&0&2&1&1&0&0&-1&-1&-2&-2\\ 2&462&14&-18&-6&0&-2&-2&6&2&2&0&0&0&-2&2&0&-2&0&0&0&-1&-1&0&0&2&2\\ 3&1540&-28&20&10&-14&4&-4&-4&0&2&2&0&0&0&0&0&-2&2&0&0&0&0&0&0&-1&-1\\ 4&4554&42&-38&0&12&-6&2&-6&-6&0&4&4&4&-2&2&0&0&0&0&0&0&0&-2&-2&0&0\\ 5&11592&-56&72&-18&0&-8&8&0&2&-2&0&0&0&0&2&-2&-2&0&0&0&2&2&0&0&0&0\\ 6&27830&86&-90&20&-16&6&-2&6&0&-4&0&-2&-2&2&0&0&0&0&2&2&0&0&-2&-2&0&0\\ 7&61686&-138&118&0&30&6&-10&-2&6&0&-2&2&2&-2&-2&-2&0&-2&2&2&0&0&2&2&0&0\\ 8&131100&188&-180&-30&0&-4&4&-12&0&2&0&-3&-3&0&0&2&2&0&-1&-1&0&0&0&0&0&0\\ 9&265650&-238&258&42&-42&-14&10&10&-10&2&6&0&0&-2&-2&0&-2&-2&0&0&2&2&0&0&0&0\\
       \bottomrule
        \multicolumn{27}{c}{}\\ \multicolumn{27}{c}{}\\ \multicolumn{27}{c}{}\\
          \toprule
  {\backslashbox{$\ell$}{$\r$}} &$ \rho_{1}$&$ \rho_{23}$&$ \rho_{45}$&$ \rho_{\overline{45}}$&$ \rho_{231}$&$ \rho_{\overline{231}}$&$ \rho_{252}$&$ \rho_{253}$&$ \rho_{483}$&$ \rho_{770}$&$ \rho_{\overline{770}}$&$ \rho_{990}$&$ \rho_{\overline{990}}$&$ \rho_{1035'}$&$ \rho_{1035}$&$ \rho_{\overline{1035}}$&$ \rho_{1265}$&$ \rho_{1771}$&$ \rho_{2024}$&$ \rho_{2277}$&$ \rho_{3312}$&$ \rho_{3520}$&$ \rho_{5313}$&$ \rho_{5544}$&$ \rho_{5796}$&$ \rho_{10395}$ \\\midrule
1&0&0&1&1&0&0&0&0&0&0&0&0&0&0&0&0&0&0&0&0&0&0&0&0&0&0\\ 
2&0&0&0&0&1&1&0&0&0&0&0&0&0&0&0&0&0&0&0&0&0&0&0&0&0&0\\ 
3&0&0&0&0&0&0&0&0&0&1&1&0&0&0&0&0&0&0&0&0&0&0&0&0&0&0\\ 
4&0&0&0&0&0&0&0&0&0&0&0&0&0&0&0&0&0&0&0&2&0&0&0&0&0&0\\ 
5&0&0&0&0&0&0&0&0&0&0&0&0&0&0&0&0&0&0&0&0&0&0&0&0&2&0\\ 
6&0&0&0&0&0&0&0&0&0&0&0&0&0&0&0&0&0&0&0&0&0&2&0&0&0&2\\ 
7&0&0&0&0&0&0&0&0&0&0&0&0&0&0&0&0&0&2&2&0&0&0&2&2&2&2\\ 
8&0&0&0&0&0&0&0&0&0&0&0&1&1&0&1&1&2&0&0&2&2&2&4&2&2&6\\ 
9&0&0&0&0&0&0&0&0&2&2&2&0&0&2&2&2&0&2&2&2&4&4&4&8&8&10\\\bottomrule
                  \end{tabular}}
     \caption{\label{Mock_decompositions} \footnotesize{The first few Fourier coefficients of the terms $q^{\frac{1}{8}+n}$ in the $q$-series $H_g(\t)$ and the corresponding representation $K_n$. Notice that the representations always come in conjugate pairs. }}
       \end{sidewaystable}
 
\newpage


\begin{thebibliography}{10}

\bibitem{sphere_packing}
J.~H. Conway and N.~J.~A. Sloane, {\em {Sphere packings, lattices and groups}}.
\newblock Springer-Verlag, 1999.

\bibitem{atlas}
J.~Conway, R.~Curtis, S.~Norton, R.~Parker, and R.~Wilson, {\em {Atlas of
  finite groups. Maximal subgroups and ordinary characters for simple groups.
  With comput. assist. from J. G. Thackray.}}
\newblock Oxford: Clarendon Press, 1985.

\bibitem{conway_norton}
J.~H. Conway and S.~P. Norton, ``{Monstrous Moonshine},'' {\em Bull. London
  Math. Soc.} {\bf 11} (1979)  308~339.

\bibitem{Mason}
G.~Mason, ``{$M_{24}$ and Certain Automorphic Forms},'' {\em Contemporary
  Mathematics} {\bf 45} (1985)  223--244.

\bibitem{DummitKisilevskyMcKay}
D.~Dummit, H.~Kisilevsky, and J.~McKay, ``Multiplicative products of
  eta-functions.,'' {\em Contemp. Math. (Finite groups - coming of age, Proc.
  CMS Conf., Montreal/Que. 1982 )} {\bf 45} (1985)  .

\bibitem{RamanunjanEta}
D.~Ford and J.~McKay, ``{Ramifications of Ramanunjan's Work on
  $\eta$-products},'' {\em Proc. Indian Acad. Sci.} {\bf 99, No.3} (1989)
  221--229.

\bibitem{generalized_moonshine}
S.~Norton, ``{Generalized Moonshine},'' {\em Proc. Symp. Pure Math} {\bf 47}
  (1987)  208--209.

\bibitem{Mason_genus0}
G.~Mason, ``{$G$-elliptic systems and the genus zero problem for $M\sb{24}$},''
  {\em Bull. Am. Math. Soc., New Ser. 25} {\bf 1} (1991)  45--53.

\bibitem{Eguchi2010}
T.~Eguchi, H.~Ooguri, and Y.~Tachikawa, ``{Notes on the $K3$ Surface and the
  Mathieu group $M_{24}$},'' \href{http://arxiv.org/abs/1004.0956}{{\tt
  1004.0956}}.

\bibitem{Eguchi2008}
T.~Eguchi and K.~Hikami, ``Superconformal algebras and mock theta
  functions,''{\em J.Phys.A} {\bf 42:304010,2009} (Dec., 2008)  ,
  \href{http://arxiv.org/abs/0812.1151}{{\tt 0812.1151}}.

\bibitem{Atish_Sameer_Don}
A.~Dabholkar, S.~Murthy, and D.~Zagier, ``{Quantum Black Hole, Wall-Crossing,
  and Mock Modular Forms},'' {\em to appear}  .

\bibitem{BriOno_MockThetaConj}
K.~Bringmann and K.~Ono, ``The {$f(q)$} mock theta function conjecture and
  partition ranks,'' \href{http://dx.doi.org/10.1007/s00222-005-0493-5}{{\em
  Invent. Math.} {\bf 165} (2006) no.~2, 243--266}.
  \url{http://dx.doi.org/10.1007/s00222-005-0493-5}.

\bibitem{Cheng2010_1}
M.~C.~N. Cheng, ``{$K3$ Surfaces, $N=4$ Dyons, and the Mathieu Group
  $M_{24}$},'' \href{http://arxiv.org/abs/1005.5415}{{\tt 1005.5415}}.

\bibitem{Gaberdiel2010}
M.~R. Gaberdiel, S.~Hohenegger, and R.~Volpato, ``{Mathieu twining characters
  for $K3$},'' \href{http://arxiv.org/abs/1006.0221}{{\tt 1006.0221}}.

\bibitem{Gaberdiel2010a}
M.~R. Gaberdiel, S.~Hohenegger, and R.~Volpato, ``{Mathieu Moonshine in the
  elliptic genus of $K3$},'' \href{http://arxiv.org/abs/1008.3778}{{\tt
  1008.3778}}.

\bibitem{Eguchi2010a}
T.~Eguchi and K.~Hikami, ``{Note on Twisted Elliptic Genus of $K3$ Surface},''
  \href{http://arxiv.org/abs/1008.4924}{{\tt 1008.4924}}.

\bibitem{Gannon_proof}
T.~Gannon, 2011.
\newblock Personal communication.

\bibitem{Cheng2011}
M.~C.~N. Cheng and J.~F.~R. Duncan, ``On rademacher sums, the largest mathieu
  group, and the holographic modularity of moonshine,''
  \href{http://arxiv.org/abs/1110.3859}{{\tt 1110.3859}}.

\bibitem{Duncan2009}
J.~F.~R. Duncan and I.~B. Frenkel, ``{Rademacher sums, moonshine and
  gravity},'' {\em To appear in Communications in Number Theory and Physics}  ,
  \href{http://arxiv.org/abs/0907.4529}{{\tt 0907.4529}}.

\bibitem{Ochanine}
S.~Ochanine, ``Sur les genres multiplicatifs définis par des intégrales
  elliptiques. (french) [on multiplicative genera defined by elliptic
  integrals],'' {\em Topology. An International Journal of Mathematics} {\bf
  26} (1987)  143.

\bibitem{Witten1987}
E.~Witten, ``{Elliptic Genera and Quantum Field Theory},''
\href{http://dx.doi.org/10.1007/BF01208956}{{\em Commun. Math. Phys.} {\bf 109}
  (1987)  525}.

\bibitem{Landweber_book}
P.~S. Landweber, ed., {\em {Elliptic curves and modular forms in algebraic
  topology}}.
\newblock Springer-Verlag, Berlin, 1988.

\bibitem{WittenInt.J.Mod.Phys.A9:4783-48001994}
E.~Witten, ``On the landau-ginzburg description of $n=2$ minimal models,'' {\em
  Int.J.Mod.Phys.A} {\bf 9:4783-4800,1994} (Int.J.Mod.Phys.A9:4783-4800,1994)
  , \href{http://arxiv.org/abs/hep-th/9304026}{{\tt hep-th/9304026}}.

\bibitem{KawaiNucl.Phys.B414:191-2121994}
T.~Kawai, Y.~Yamada, and S.-K. Yang, ``Elliptic genera and n=2 superconformal
  field theory,'' {\em Nucl.Phys.B} {\bf 414:191-212,1994}
  (Nucl.Phys.B414:191-212,1994)  ,
  \href{http://arxiv.org/abs/hep-th/9306096}{{\tt hep-th/9306096}}.

\bibitem{Gritsenko2007}
V.~Gritsenko, ``{Elliptic genus of Calabi-Yau manifolds and Jacobi and Siegel
  modular forms},'' \href{http://arxiv.org/abs/math/9906190}{{\tt
  math/9906190}}.

\bibitem{Eguchi1989}
T.~Eguchi, H.~Ooguri, A.~Taormina, and S.-K. Yang, ``{Superconformal Algebras
  and String Compactification on Manifolds with SU(N) Holonomy},''
\href{http://dx.doi.org/10.1016/0550-3213(89)90454-9}{{\em Nucl. Phys.} {\bf
  B315} (1989)  193}.

\bibitem{Eguchi1987}
T.~Eguchi and A.~Taormina, ``{Unitary Representations of $N=4$ Superconformal
  Algebra},''
\href{http://dx.doi.org/10.1016/0370-2693(87)91679-0}{{\em Phys. Lett.} {\bf
  B196} (1987)  75}.

\bibitem{Eguchi1988}
T.~Eguchi and A.~Taormina, ``{Character Formulas for the $N=4$ Superconformal
  Algebra},''
\href{http://dx.doi.org/10.1016/0370-2693(88)90778-2}{{\em Phys. Lett.} {\bf
  B200} (1988)  315}.

\bibitem{Eguchi1988a}
T.~Eguchi and A.~Taormina, ``{On the Unitary Representations of $N=2$ and $N=4$
  Superconformal Algebras},''
\href{http://dx.doi.org/10.1016/0370-2693(88)90360-7}{{\em Phys. Lett.} {\bf
  B210} (1988)  125}.

\bibitem{Eguchi2009a}
T.~Eguchi and K.~Hikami, ``{Superconformal Algebras and Mock Theta Functions 2.
  Rademacher Expansion for $K3$ Surface},''{\em Communications in Number Theory
  and Physics} {\bf 3,} (Apr., 2009)  531--554,
  \href{http://arxiv.org/abs/0904.0911}{{\tt 0904.0911}}.

\bibitem{zwegers}
S.~Zwegers, {\em {Mock Theta Functions}}.
\newblock PhD thesis, Utrecht University, 2002.

\bibitem{Mukai}
S.~Mukai, ``{Finite groups of automorphisms of $K3$ surfaces and the Mathieu
  group},'' {\em Invent. Math.} {\bf 94, No.1} (1988)  183--221.

\bibitem{Kondo}
S.~Kondo, Shigeyuki~(Mukai, ``{Niemeier lattices, Mathieu groups, and finite
  groups of symplectic automorphisms of $K3$ surfaces. With an appendix by
  Shigeru Mukai.},'' {\em Duke Math. J} {\bf 92, No.3} (1998)  593--603.

\bibitem{Taormina2010}
A.~Taormina and K.~Wendland, ``{The symmetries of the tetrahedral Kummer
  surface in the Mathieu group $M_{24}$},''
  \href{http://arxiv.org/abs/1008.0954}{{\tt 1008.0954}}.

\bibitem{Gaberdiel2011}
M.~R. Gaberdiel, S.~Hohenegger, and R.~Volpato, ``{Symmetries of $K3$ sigma
  models},'' \href{http://arxiv.org/abs/1106.4315}{{\tt 1106.4315}}.

\bibitem{DavidJHEP0606:0642006}
J.~R. David, D.~P. Jatkar, and A.~Sen, ``{Product Representation of Dyon
  Partition Function in CHL Models},'' {\em JHEP} {\bf 0606:064,2006} (JHEP
  0606:064,2006)  , \href{http://arxiv.org/abs/hep-th/0602254}{{\tt
  hep-th/0602254}}.

\bibitem{Govindarajan2008}
S.~Govindarajan and K.~G. Krishna, ``{Generalized Kac-Moody Algebras from CHL
  dyons},''{\em JHEP} {\bf 0904:032,2009} (July, 2008)  ,
  \href{http://arxiv.org/abs/0807.4451}{{\tt 0807.4451}}.

\bibitem{Sen2010}
A.~Sen, ``{Discrete Information from CHL Black Holes},''
  \href{http://arxiv.org/abs/1002.3857}{{\tt 1002.3857}}.

\bibitem{DVV}
R.~Dijkgraaf, E.~Verlinde, and H.~Verlinde, ``{Counting Dyons in N=4 String
  Theory},'' {\em Nucl.Phys.B} {\bf 484:543-561,1997}
  (Nucl.Phys.B484:543-561,1997)  ,
  \href{http://arxiv.org/abs/hep-th/9607026}{{\tt hep-th/9607026}}.

\bibitem{Dijkgraaf1997}
R.~Dijkgraaf, G.~W. Moore, E.~P. Verlinde, and H.~L. Verlinde, ``{Elliptic
  genera of symmetric products and second quantized strings},''
  \href{http://dx.doi.org/10.1007/s002200050087}{{\em Commun. Math. Phys.} {\bf
  185} (1997)  197--209},
\href{http://arxiv.org/abs/hep-th/9608096}{{\tt arXiv:hep-th/9608096}}.

\bibitem{igusa}
J.-I. Igusa, ``{On Siegel modular forms of genus two (I)},'' {\em Am. J. Math.}
  {\bf 84} (1962)  175--200.

\bibitem{igusa2}
J.-I. Igusa, ``{On Siegel modular forms of genus two (II)},'' {\em Am. J.
  Math.} {\bf 86} (1964)  392--412.

\bibitem{DavidJHEP0611:0722006}
J.~R. David and A.~Sen, ``{CHL Dyons and Statistical Entropy Function from
  D1-D5 System},'' {\em JHEP} {\bf 0611:072,2006} (JHEP 0611:072,2006)  ,
  \href{http://arxiv.org/abs/hep-th/0605210}{{\tt hep-th/0605210}}.

\bibitem{ShihJHEP0610:0872006}
D.~Shih, A.~Strominger, and X.~Yin, ``{Recounting Dyons in N=4 String
  Theory},'' {\em JHEP} {\bf 0610:087,2006} (JHEP 0610:087,2006)  ,
  \href{http://arxiv.org/abs/hep-th/0505094}{{\tt hep-th/0505094}}.

\bibitem{Sen2007c}
A.~Sen, ``Black hole entropy function, attractors and precision counting of
  microstates,''{\em Gen.Rel.Grav.} {\bf 40:2249-2431,2008} (Aug., 2007)  ,
  \href{http://arxiv.org/abs/0708.1270}{{\tt 0708.1270}}.

\bibitem{Cheng2008b}
M.~C.~N. Cheng, ``The spectra of supersymmetric states in string theory,''
  \href{http://arxiv.org/abs/0807.3099}{{\tt 0807.3099}}.

\bibitem{SenJHEP0705:0392007}
A.~Sen, ``{Walls of Marginal Stability and Dyon Spectrum in N=4 Supersymmetric
  String Theories},'' {\em JHEP} {\bf 0705:039,2007} (JHEP 0705:039,2007)  ,
  \href{http://arxiv.org/abs/hep-th/0702141}{{\tt hep-th/0702141}}.

\bibitem{Cheng2007a}
M.~C.~N. Cheng and E.~Verlinde, ``{Dying Dyons Don't Count},''
  \href{http://dx.doi.org/10.1088/1126-6708/2007/09/070}{{\em JHEP} {\bf 09}
  (2007)  070},
\href{http://arxiv.org/abs/0706.2363}{{\tt arXiv:0706.2363 [hep-th]}}.

\bibitem{Dabholkar2008}
A.~Dabholkar, D.~Gaiotto, and S.~Nampuri, ``{Comments on the spectrum of CHL
  dyons},'' \href{http://dx.doi.org/10.1088/1126-6708/2008/01/023}{{\em JHEP}
  {\bf 01} (2008)  023},
\href{http://arxiv.org/abs/hep-th/0702150}{{\tt arXiv:hep-th/0702150}}.

\bibitem{borcherds_singular}
R.~E. Borcherds, ``{Automorphic forms with singularities on Grassmannians},''.

\bibitem{Grit_Nik}
V.~A. Gritsenko and V.~V. Nikulin, ``{Siegel automorphic form corrections of
  some Lorentzian Kac-Moody algebras},'' {\em Amer.J.Math.} {\bf 119} (1991)
  181.

\bibitem{Gritsenko2008b}
V.~A. Gritsenko and V.~V. Nikulin, ``Automorphic forms and lorentzian kac-moody
  algebras. part ii,'' \href{http://arxiv.org/abs/alg-geom/9611028}{{\tt
  alg-geom/9611028}}.

\bibitem{Carnahan2008}
S.~Carnahan, ``Generalized moonshine i: Genus zero functions,''{\em Algebra and
  Number Theory} {\bf 4:6} (Dec., 2008)  649--679,
  \href{http://arxiv.org/abs/0812.3440}{{\tt 0812.3440}}.

\bibitem{ganter1}
N.~Ganter, ``{Hecke operators in equivariant elliptic cohomology and
  generalized moonshine},'' {\em [Harnad, John (ed.) et al., Groups and
  symmetries. From Neolithic Scots to John McKay. Selected papers of the
  conference, Montreal, Canada, April 27--29, 2007. Providence, RI: American
  Mathematical Society (AMS). CRM Proceedings and Lecture Notes 47, 173-209
  (2009; Zbl 1206.55006)]}  .

\bibitem{DavidJHEP0701:0162007}
J.~R. David, D.~P. Jatkar, and A.~Sen, ``{Dyon Spectrum in Generic N=4
  Supersymmetric $Z_N$ Orbifolds},'' {\em JHEP} {\bf 0701:016,2007} (JHEP
  0701:016,2007)  , \href{http://arxiv.org/abs/hep-th/0609109}{{\tt
  hep-th/0609109}}.

\bibitem{DabholkarJHEP0711:0772007}
A.~Dabholkar and S.~Nampuri, ``{Spectrum of Dyons and Black Holes in CHL
  orbifolds using Borcherds Lift},'' {\em JHEP} {\bf 0711:077,2007}
  (JHEP0711:077,2007)  , \href{http://arxiv.org/abs/hep-th/0603066}{{\tt
  hep-th/0603066}}.

\bibitem{Sen2009}
A.~Sen, ``{A Twist in the Dyon Partition Function},''
  \href{http://arxiv.org/abs/0911.1563}{{\tt 0911.1563}}.

\bibitem{Govindarajan2009}
S.~Govindarajan and K.~G. Krishna, ``{BKM Lie superalgebras from dyon spectra
  in $Z_N$-CHL orbifolds for composite N},''
  \href{http://arxiv.org/abs/0907.1410}{{\tt 0907.1410}}.

\bibitem{Govindarajan2011}
S.~Govindarajan, ``Unravelling mathieu moonshine,''
  \href{http://arxiv.org/abs/1106.5715}{{\tt 1106.5715}}.

\bibitem{Eguchi2011}
T.~Eguchi and K.~Hikami, ``{Twisted elliptic genus for $K3$ and Borcherds
  product},'' \href{http://arxiv.org/abs/1112.5928}{{\tt 1112.5928}}.

\bibitem{Gritsenko2008}
V.~Gritsenko and F.~Clery, ``The siegel modular forms of genus 2 with the
  simplest divisor,'' \href{http://arxiv.org/abs/0812.3962}{{\tt 0812.3962}}.

\bibitem{JatkarJHEP0604:0182006}
D.~P. Jatkar and A.~Sen, ``{Dyon Spectrum in CHL Models},'' {\em JHEP} {\bf
  0604:018,2006} (JHEP 0604:018,2006)  ,
  \href{http://arxiv.org/abs/hep-th/0510147}{{\tt hep-th/0510147}}.

\bibitem{Tuite2001}
M.~P. Tuite, ``Genus two meromorphic conformal field theory,'' {\em CRM
  Proceedings \& Lecture Notes} {\bf 30,} (2001)  231--251.,
  \href{http://arxiv.org/abs/math/9910136}{{\tt math/9910136}}.

\bibitem{DabholkarJHEP0712:0872007}
A.~Dabholkar and D.~Gaiotto, ``{Spectrum of CHL Dyons from Genus-Two Partition
  Function},'' {\em JHEP} {\bf 0712:087,2007} (JHEP0712:087,2007)  ,
  \href{http://arxiv.org/abs/hep-th/0612011}{{\tt hep-th/0612011}}.

\bibitem{Govindarajan2010a}
S.~Govindarajan, ``Brewing moonshine for mathieu,''
  \href{http://arxiv.org/abs/1012.5732}{{\tt 1012.5732}}.

\bibitem{feingold_frenkel}
A.~J. Feingold and I.~B. Frenkel, ``{A hyperbolic Kac-Moody algebra and the
  theory of Siegel modular forms of genus 2},'' {\em J. Math. Ann.} {\bf 263}
  (1983)  87--144.

\bibitem{feingold_nicolai}
A.~J. Feingold and H.~Nicolai, ``Subalgebras of hyperbolic {K}ac-{M}oody
  algebras,'' in {\em Kac-{M}oody {L}ie algebras and related topics}, vol.~343
  of {\em Contemp. Math.}, pp.~97--114.
\newblock Amer. Math. Soc., Providence, RI, 2004.

\bibitem{GriNik_IguLKM}
V.~A. Gritsenko and V.~V. Nikulin, ``Igusa modular forms and ``the simplest''
  {L}orentzian {K}ac-{M}oody algebras,''
  \href{http://dx.doi.org/10.1070/SM1996v187n11ABEH000171}{{\em Mat. Sb.} {\bf
  187} (1996) no.~11, 27--66}.
  \url{http://dx.doi.org/10.1070/SM1996v187n11ABEH000171}.

\bibitem{eichler_zagier}
M.~Eichler and D.~Zagier, {\em {The theory of Jacobi forms}}.
\newblock Birkh\"auser, 1985.

\bibitem{modi}
\url{http://modi.countnumber.de/index.php?chap=ell.newforms/ell.newforms.html}.

\bibitem{from_number}
C.~Itzykson, ed., {\em {From Number Theory to Physics}}.
\newblock Springer, 1992.

\bibitem{aoki}
H.~Aoki and T.~Ibukiyama, ``{Simple graded rings of Siegel modular forms,
  differential operators and Borcherds products},'' {\em Int. J. Math.} {\bf
  16, No. 3} (2005)  249--279.

\end{thebibliography}

\providecommand{\href}[2]{#2}\begingroup\raggedright\endgroup

\end{document}